\documentclass[leqno,12pt]{article}

\usepackage{hyperref} 
\usepackage{amsthm}

\usepackage{graphicx}
\usepackage{amssymb}
\usepackage{amsfonts}
\usepackage{latexsym}
\usepackage{amscd}
\usepackage[mathscr]{euscript}
\usepackage{xy} \xyoption{all}
\usepackage[utf8]{inputenc} 
\usepackage{comment}
\usepackage{caption}
\usepackage{amsmath}
\usepackage{mathrsfs}
\usepackage{mathtools}

\makeatletter
\newcommand{\leqnomode}{\tagsleft@true}
\newcommand{\reqnomode}{\tagsleft@false}
\makeatother

\usepackage{tikz}
\usetikzlibrary{shapes,arrows}
\usepackage[utf8]{inputenc}
\usepackage{color}
\newcommand{\DistTo}{\xrightarrow{
   \,\smash{\raisebox{-0.65ex}{\ensuremath{\scriptstyle\sim}}}\,}}

\usepackage{setspace}

\def\Q{\mathbb Q}

\def\Z{\mathbb Z}

\theoremstyle{plain}
\newtheorem{theorem}{Theorem}[section]
\newtheorem{proposition}[theorem]{Proposition}
\newtheorem{lemma}[theorem]{Lemma}
\newtheorem{corollary}[theorem]{Corollary}

\theoremstyle{definition}

\newtheorem{example}[theorem]{Example}
\newtheorem{remark}[theorem]{Remark}
\newcommand{\Mod}[1]{\ (\mathrm{mod}\ #1)}
\def\Q{\mathbb Q}

\def\Z{\mathbb Z}

\usepackage{calc}

\usepackage{accents}

 \newcommand{\GL}{\operatorname{GL}}
 \newcommand{\SL}{\operatorname{SL}}
 \newcommand{\Div}{\operatorname{Div}}
  
 \newcommand{\Aut}{\operatorname{Aut}}
 \newcommand{\NS}{\operatorname{NS}}
 \newcommand{\rank}{\operatorname{rank}}
 \newcommand{\Hom}{\operatorname{Hom}}
 \newcommand{\End}{\operatorname{End}}
 \newcommand{\charec}{\operatorname{char}}
 \newcommand{\gen}{\operatorname{gen}}
   \newcommand{\disc}{\operatorname{disc}}
    
     \newcommand{\cont}{\operatorname{cont}}
      \newcommand{\odd}{\operatorname{odd}}
     \newcommand{\ev}{\operatorname{ev}}
     \newcommand{\adj}{\operatorname{adj}}
     \newcommand{\cl}{\operatorname{cl}}
     
\newcommand{\sps}{\vspace{3pt}}   
\newcommand{\spm}{\vspace{6pt}}

\newenvironment{romanlist}
{\begin{enumerate}[label=(\roman*), leftmargin=*]}
{\end{enumerate}}

\title{The classification of the refined Humbert invariant for curves of genus $2$
}

\author{Harun K{\i}r\footnote{Current address: ENS de Lyon, UMPA, UMR 5669, Lyon, France. \href{mailto:harun.kir@ens-lyon.fr}{harun.kir@ens-lyon.fr}}}

\begin{document}


\maketitle


\begin{abstract}
The refined Humbert invariant is a positive definite quadratic form intrinsically attached to a curve $C$ of genus 2. This invariant is an algebraic generalization of the (usual) Humbert invariant. This invariant is useful because many geometric properties of $C$ are reflected in the arithmetic properties of this invariant. The purpose of this paper is to complete the classification of this invariant when the Jacobian $J_C$ of $C$ is isogenous to a product of an elliptic curve with complex multiplication.
\end{abstract}

\section{Introduction}



If $C/K$ is a (smooth) curve of genus 2 over an algebraically closed field $K$ of characteristic $0$, then it comes equipped with a canonical quadratic form $q_C$ called its \textit{refined Humbert invariant}; cf.\ \cite{kani1994elliptic}, \cite{kani2014jacobians} or \S\ref{refinedhumbert} below. It is introduced by  Kani \cite{kani1994elliptic} as an algebraic generalization of the (usual) Humbert invariant (see Section 5 of \cite{kani1994elliptic}).   It encodes many special geometric properties of $C$, and  thus, several geometric
properties of $C$ are translated into arithmetic properties of $q_{C}$. This makes it worthwhile to look at this invariant. An elegant illustration of making use of the refined Humbert invariant can be found in Kani \cite{kani1994elliptic}, \cite{kani2014jacobians}, \cite{MJ}, \cite{SubcoversofCurves}.

It is interesting and useful to classify the quadratic forms $q$, equivalent to a refined Humbert invariant $q_{C}$, for some curve $C$ of genus 2.  Kani provided such a classification when the Jacobian $J_C$ of $C$ is isogenous to a product $E\times E$ of an elliptic curve $E$ which does not have complex multiplication (CM) (i.e., $\End(E) = \Z$); see Theorems 1 and 3 of \cite{ESCII}.  Recently, Kani \cite{refhum}  gave a similar classification for the cases in which $J_C \sim E\times E$ are isogenous, where $E$ is a CM elliptic curve and $q_{C}$ is a \textit{primitive} form. 
   
   The main aim of this article is to give a  classification when  $J_C \sim E \times E$, where $E$ is a CM elliptic curve and $q_{C}$ is \textit{not} a primitive form, i.e., \textit{imprimitive} form. Then, the classification will be complete when $J_C \sim E \times E$, with CM elliptic curve $E$.  To this end, we prove two main interconnected structural theorems. The first main result of this article is as follows. 

\begin{theorem}
\label{theorem1}
Assume that an imprimitive integral ternary quadratic form $f_1$  is equivalent to a refined Humbert invariant $q_C$,  for some curve $C/K$ of genus $2$. If a form $f_2$ is \emph{genus-equivalent} to $f_1,$ then there exists a curve $C'/K$ of genus $2$ on the Jacobian $J_C$ of $C$ such that $f_2$  is equivalent to the refined Humbert invariant $q_{C'}$.
\end{theorem}

This theorem together with that of \cite{refhum} can be applied to an interesting question, namely,  to giving a formula for the number of isomorphism classes of (smooth) curves of genus 2 lying on an abelian surface $A$ in many cases. This work is available in Kani's personal webpage.

Moreover, Theorem \ref{theorem1} constitutes a key tool to prove the second main theorem (see Theorem \ref{maintheorem} below), which gives a classification of the imprimitive (integral) ternary quadratic forms $f$  which are equivalent to a refined Humbert invariant $q_{C},$ for some curve $C$ of genus 2.  For this, let us consider positive definite ternary quadratic forms $f(x,y,z)$ satisfying the following  two conditions:
\begin{align}
\label{classificationconditions1}
  \frac{1}{2}f& \text{ is an improperly primitive form; (see \S\ref{ternaryforms} below)} \\
 f(x_0,y_0,z_0)& \;=\; (2n)^2 \textrm{ for some } x_0,y_0,z_0,n\in\Z \textrm{ with } \gcd(n,\disc(f))\;=\;1, \label{classificationconditions2}
\end{align} where $\disc(f)$ is defined as in \cite[p.~2]{watson1960integral}. Then we have the following existence theorem.

\begin{theorem}
\label{maintheorem}
Let $f$ be an imprimitive positive definite integral ternary quadratic form. Then there is a curve $C/K$ of genus $2$ such that 
 $f$ is equivalent to a refined Humbert invariant $q_C$  if and only if
$f$ satisfies  conditions \emph{(\ref{classificationconditions1})} and \emph{(\ref{classificationconditions2})}. If this is the case, then the Jacobian $J_C$ of $C$ is isogenous to a product $E \times E$ of an elliptic curve $E/K$ with complex multiplication.
\end{theorem}

 Theorem \ref{maintheorem}  has interesting applications; cf.\ Section \ref{section: applications} below. To mention some of them,  recall first from \cite{kani1994elliptic}, \cite{MJ} that the refined Humbert invariant $q_{(A,\theta)}$ is defined for any principally polarized abelian surface $(A,\theta) \in \mathcal{A}_2(K)$ (cf.\ \S\ref{refinedhumbert} below), and  we have by definition that $q_C:=q_{(J_C,\theta_C)}$, if $C/K$ is a curve of genus 2. Recall then that Kani \cite{MJ} introduced the notion of a \textit{generalized Humbert scheme} $H(q) \subset \mathcal{A}_2(K)$ that is associated to a given quadratic form $q$. This scheme is useful to determine how the curves of genus $2$ on the product abelian surfaces are distributed in the moduli space $\mathcal{M}_2(K)$ of genus $2$ curves; see \cite{MJ}. This scheme is defined  using the refined Humbert invariant, and it is a generalization of a (usual) \textit{Humbert surface $H_N$ of invariant} $N$; see \S3 of \cite{MJ} or \S \ref{section: applications} below (also, see van der Geer \cite{van2012hilbert}, Ch. IX, for a contemporary treatment of Humbert surfaces).

 To understand the components of the intersection $H_{N^2}\cap H_{M^2}$ of two Humbert surfaces with square invariant, Kani (see Question 6 of \cite{SubcoversofCurves}) raised the following question. For which integral binary quadratic forms $q$ is $H(q) \neq \varnothing$ provided $H(q) \subset H_{N^2}$, for some $N$? He provided a  solution to this problem; see Theorem 8.1 of \cite{SubcoversofCurves}. As a generalization of this question, we ask for which integral binary quadratic forms $q$ is $H(q) \neq \varnothing$ ?   One can obtain a partial solution when a binary form does not represent a square and its discriminant is square-free by using Theorem 5.2 of \cite{Hashimoto}. We give a complete answer to this question  for all cases with a different method that constructs \textit{CM points} on $H(q)$.

\begin{corollary}
\label{cor: binaryform_H(q)}
    Let $q$ be a positive definite integral binary quadratic form. Then we have that
    $$
    H(q) \;\neq\; \varnothing \quad \Leftrightarrow \quad q \;\equiv\; 0,1 \Mod{4}.
    $$
\end{corollary}

Theorem \ref{maintheorem} and Theorem 1 of \cite{refhum} play a key role in the proof of this corollary. The properties of the set $H(q)$ were analyzed for binary quadratic forms $q$ that represent a perfect square in \cite{MJ}, \cite{SubcoversofCurves}. 
 Thus, this corollary leads to numerous questions regarding the structure of $H(q)$ in cases where $q$ is a positive definite binary quadratic form that does not represent a square. (In this case, $H(q)$ has a close relation with the Shimura curves as will be shown in the author's PhD thesis.)  For example, one of the interesting questions is to determine the irreducible components of $H(q)$, for a binary quadratic form $q$ because this can be used to describe the components of the intersection of two Humbert surfaces; cf.\ Remark 8 of \cite{MJ} and/or Question 5 (and its partial answer) of \cite{SubcoversofCurves}. Note that the problem of studying the components of this intersection was briefly addressed on \cite[p.~214]{van2012hilbert}. Thus, Corollary \ref{cor: binaryform_H(q)} can be seen as a contribution to this discussion in the view of the nice formula/algorithm of the intersection of two Humbert surfaces in terms of the unions of the $H(q)$, for binary forms $q$; cf.\ Proposition 5.2.d of \cite{SubcoversofCurves}.  

Related to the intersection of two Humbert surfaces,  Corollary \ref{cor: binaryform_H(q)} has an immediate application, which can be seen as a contribution of  McMullen \cite{Mcmullen}'s question. More precisely, he posed a question regarding the description of the intersection of Humbert surfaces $H_N \cap H_M$.  He responded to this problem for the intersection $H_N \cap H_1$. Concerning this question, the first natural task might be to determine whether the intersection  $H_N \cap H_M$ is empty for any (positive) $N, M$. Franciosi, Pardini and Rollenske \cite{Pardini_Marco_Rollenske} showed that $H_{N^2} \cap H_{M^2}\neq \varnothing$. By a different method, when $H_N \neq \varnothing$, Kani showed that $H_N \cap H_{M^2}\neq \varnothing$ in Corollary 8.2 of \cite{SubcoversofCurves}. From Corollary \ref{cor: binaryform_H(q)}, we can deduce the following result.
\begin{corollary}
\label{cor: intersection_humbert_surfaces}
   Let $H_N$ be the Humbert surface of invariant $N$. For non-empty Humbert surfaces $H_N$ and $H_M$ we have that the intersection $H_N \cap H_M \neq \varnothing$. 
    
    In other words, we have that if $N, M >0$ and $N, M \equiv 0, 1 \Mod{4}$, then  $H_N \cap H_M \neq \varnothing$.
\end{corollary}

For a further discussion, it might be interesting to mention an illustration of Theorem \ref{maintheorem} for the intersection of Humbert surfaces as follows:
\begin{proposition}
\label{p: infinite_intersection_surfaces}
    Let $C$ denote the curve of genus $2$ defined by $y^2 = 9x^5 + 40x^3 + 45x$. Then we have that
\begin{equation*}
\label{eq: infinite_intersection_surfaces}
(J_C, \theta_C) \;\in\; \bigcap_N \;\; H_{4N} ,\ \text{ where } N \text{ runs over all positive odd squarefree integers}.
\end{equation*}
\end{proposition}
This result uses a special ternary quadratic form $q$ which represents all odd integers in Kaplansky's list \cite{kaplansky} and also uses the curves listed in Table 4 of \cite{gelin2019principally}.

Our main classification theorem (Theorem \ref{maintheorem}) leads to other applications beyond the scope of the present article.  Let us mention some of these applications here. Firstly, the curves of genus 2 in terms of the refined Humbert invariants can be listed based on their automorphism groups. This was studied in the author's forthcoming paper, and it is available in arXiv:2310.19076. Secondly, with a joint work of E. Kani,  for a given ternary form $q$, we provided a formula for the number of the isomorphism classes of curves  $C$ of genus 2 such that $q_{C}$ is equivalent 
to $q$ in another forthcoming 
paper, which is available in author's personal webpage. Finally, we characterize all \textit{CM points} on the \textit{Shimura curves} at the intersections of two Humbert surfaces considered in the papers of Runge \cite{Runge}, Hashimoto \cite{Hashimoto}, Lin and Yang \cite{quaternionic} (and other papers); cf.\ Remark \ref{rem: ref_hum_inv_relations_literature} below. This topic will appear in the PhD thesis of the author.

Although proving the necessary condition of Theorem \ref{maintheorem} is easy, the proof of the sufficient condition is quite long and occupies most of the present article. Our strategy has some similarities to that of \cite{refhum}. Therefore, throughout our work, the same notations and conventions are mostly used as in this article. 

The organization of the paper is as follows. In Section \ref{ternaryforms}, we review some fundamental concepts and terminology in \cite{brandt1952mass}, \cite{brandt1951zahlentheorie}, \cite{dicksonsbook}, \cite{jones1950arithmetic}, \cite{smith} and \cite{watson1960integral} regarding the theory of integral ternary quadratic forms. 

In Section \ref{refinedhumbert}, we define as in \cite{kani1994elliptic}, \cite{MJ} the main ingredient of our paper, namely the refined Humbert invariant, and recall its basic properties; cf.\ Propositions \ref{imprimitiveform} and  \ref{Q2aoddcase}. We also give the relation between the refined Humbert invariant (when it is ternary) and the binary forms of the type studied in \cite{MJ}; cf.\ Proposition 29 of \cite{MJ}. This relation is important for the proof of our main theorems.


 We first prove the necessity part of Theorem \ref{maintheorem}  in Section \ref{section: quadratic forms in Q(A)}; cf.\ Propositions \ref{invariants} and \ref{necessaryconditionforii}. Then we canonically give a construction of exact principally polarized abelian surfaces. This has many useful results, for example, it gives that all ternary forms satisfying (\ref{classificationconditions1}) and (\ref{classificationconditions2})  lie in a single genus in most cases; cf.\ Theorem \ref{cor: Theta_lies_genus}. This is the main result of Section \ref{section: quadratic forms in Q(A)}, and it requires to study of the properties of the imprimitive ternary quadratic forms equivalent to a refined Humbert invariant.

To give the proof of Theorem \ref{theorem1} in Section \ref{section: proofs_main results}, we aim to apply Theorem 34 of \cite{dicksonsbook}, and thus, our purpose is to show that if an imprimitive ternary form $q_{(A,\theta)}$ is genus equivalent to a form $f$, then they satisfy the hypothesis stated in this theorem. This requires constructing a  binary quadratic form with a suitable discriminant properly represented by the forms $f$ and $q_{(A,\theta)}$. This construction involves using Proposition \ref{frepresentsphi} and Theorem \ref{phiprincipalgenus}. Then Theorem \ref{maintheorem} follows from Theorem \ref{theorem1}.

In the last section, we prove Corollaries \ref{cor: binaryform_H(q)}, \ref{cor: intersection_humbert_surfaces} and Proposition \ref{p: infinite_intersection_surfaces}, and give explicit examples to illustrate Theorem \ref{maintheorem}.
\section{Ternary quadratic forms} 
\label{ternaryforms}

 In this article, we mostly work on positive definite binary and ternary integral quadratic forms as was mentioned in the introduction. For this, we want to use  results of Brandt \cite{brandt1952mass}, \cite{brandt1951zahlentheorie}, Cox \cite{davidcox}, Dickson \cite{dicksonsbook}, Jones \cite{jones1950arithmetic}, Smith \cite{smith} and Watson \cite{watson1960integral}. It is necessary for our aim to identify the assigned characters of a ternary quadratic form and its reciprocal and to calculate their values, and so we briefly mention some notations and concepts in the sense of those authors.

 An \textit{integral quadratic form} in $n$ variables is a polynomial $f(x_1,\ldots, x_n) $ of the form 
 \begin{equation*}
\label{quadraticform}
    f \;=\; f(x_1,\ldots, x_n) \;=\; \sum_{1\leq i\leq j\leq n}a_{ij}x_ix_j, \textrm{ where } a_{ij}\in \Z
\end{equation*} 
 as defined in Watson \cite{watson1960integral} (and in Brandt \cite{brandt1951zahlentheorie}). 
Note that an integral quadratic form $f$ (in the above sense) may not be an integral form in the sense of Dickson \cite{dicksonsbook} and Smith \cite{smith} because they additionally assume that the \textit{non-diagonal coefficients} of $f$ are even, i.e., $2\mid a_{ij},$ for all $i<j.$

The \textit{content} of the form $f$ is the greatest common divisor of its coefficients $a_{ij}$, and is denoted by $\cont(f).$ If $\cont(f)=1$, we say that it is a \textit{primitive} form; otherwise, it is called \textit{imprimitive}. For the older terminology in \cite{dicksonsbook}, \cite{smith}, $f$ is called \textit{properly primitive} if it is primitive (in the above sense) and if $2\mid a_{ij},$ for all $i<j.$ Also, $f$ is called 
\textit{improperly primitive} if $\frac{1}{2}f$ is primitive (in the above sense) and if there exists some pair $(i,j)$ with $i<j$ such that $\frac{1}{2}a_{ij}$ is odd.

By $\disc(f)$ and $\det(f),$ we mean the \textit{discriminant} and the \textit{determinant} of an integral quadratic form $f$ respectively as defined in \cite[p.~2]{watson1960integral}. Additionally, if we let $\det^{\text{D}}(f)$ the determinant of a positive definite integral binary form $f(x,y)=ax^2+2txy+by^2$ in \textit{Dickson's sense} (cf.\ \cite[p.~3]{dicksonsbook}),  then we have that 
\begin{equation}
\label{dicksondeterminant}
    {\det}^{\text{D}}(f) \;=\; \frac{-\det(f)}{4} \;=\; t^2-ab.
\end{equation}

Let $f_1$ and $f_2$ be two integral quadratic forms. If $f_1\sim_pf_2$ for all primes $p$ including $p=\infty,$ then we say that $f_1$ and $f_2$ are 
\textit{genus-equivalent} (or \textit{semi-equivalent}); cf.\ \cite[p.~106]{jones1950arithmetic}. Note that this is equivalent to the definition which was given in Watson \cite[p.~72]{watson1960integral}; cf.\ Theorem 40 of \cite{jones1950arithmetic} (also see Theorem 50 of \cite{watson1960integral}). In this case, we say that they are in the same \textit{genus}. Let $\gen(f)$ denote the set of equivalence classes of integral quadratic forms $f'$ which are genus-equivalent to a given integral quadratic form $f.$

Note that if $f_1$ and $f_2$ are positive definite improperly primitive integral ternary forms, then they are genus-equivalent (in the above sense) if and only if they are genus-equivalent in Smith's sense \cite{smith} by Theorem 40 of \cite{jones1950arithmetic} and \S12 of \cite{smith}. Therefore, we study certain techniques for deciding when two positive definite integral ternary forms $f_1$ and $f_2$ are genus-equivalent by using the method of Smith \cite{smith} in order to prove our main results. This method was improved in \cite{brandt1952mass}, \cite{brandt1951zahlentheorie} by shortening the case distinctions. To avoid many case distinctions, we make use of this improvement. For this aim, we introduce some terminology for quadratic forms here.

Throughout the article, we usually talk about the integral quadratic forms, and we often drop the word integral when it is clear in the context.

If $f$ is properly or improperly primitive, then it has a \textit{reciprocal} quadratic form $F_f$ in Smith's sense as defined on \cite[p.~455]{smith} (or on \cite[p.~7]{dicksonsbook}). If $f$ is a primitive form, then it has \textit{reciprocal} as defined by Brandt \cite[p.~336]{brandt1951zahlentheorie}; this is denoted by $F_f^B.$ We will discuss the relations of these reciprocals when $f$ is an improperly ternary form; see Proposition \ref{F=F^B} below.

To decide when two positive definite ternary forms are genus-equivalent by using Smith's method, we determine the \textit{assigned characters} of $f_1$ and $F_{f_1}$ and those of $f_2$ and $F_{f_2},$ and then we compare the values of the assigned characters. In particular, we use Smith's table \cite{smith} to decide the values of the \textit{supplementary assigned characters}.

Let us recall the notations of basic invariants as Smith and Brandt defined, and discuss how they are related to each other. Let the \textit{basic (genus) invariants} $\Omega$ and $\Delta$ be as in \cite[p.~455]{smith} (and also see \cite[p.~7]{dicksonsbook}), and let the \textit{basic (genus) invariants} $I_1$ and $I_2$ be as in Brandt \cite[p.~316]{brandt1952mass}, \cite[p.~336]{brandt1951zahlentheorie}. We first derive from \cite{brandt1952mass} their relations here, and later,  we calculate them for the forms we studied on; cf.\ Proposition \ref{invariants} below. 

\begin{proposition} 
\label{F=F^B}
Let $f$ be a positive definite improperly primitive ternary form. Then $F_f$ is properly primitive, and $F_f=F_{f/2}^B.$ Moreover, we have that
\begin{align}
\label{I1odd}
     I_1(f/2)& \   \text{ is negative odd and } \quad 16\mid I_2(f/2),\\
      \label{discf}
         \disc(f/2)& \;=\; \frac{I_1(f/2)^2I_2(f/2)}{16}, \quad \text{ and so,}\\
         \label{I1Omega}
I_1(f/2)& \;=\; -\Omega_f \quad \text{ and } \quad I_2(f/2) \;=\; -8\Delta_f.
\end{align}
\end{proposition}
\begin{proof}
Since $f$ is an improperly primitive ternary form, the primitive form $f/2$ satisfies the conditions of Case I of \cite[p.~316]{brandt1952mass}, and so $F_{f/2}^B=F_{f}$ is properly primitive and $I_1(f/2)$ is negative odd and $16\mid I_2(f/2)$ from the relations in \cite{brandt1952mass},  loc. cit. By the identity on \cite[p.~316]{brandt1952mass}, Eq.~(\ref{discf}) holds. Moreover, we get from the relations of Type I of \cite{brandt1952mass}, loc. cit., that $I_1(f/2)=-\Omega_f \text{ and } I_2(f/2)=-8\Delta_f.$ Hence, all the assertions have been verified.
\end{proof}



Let us put $ \chi_{\ell}(x) = \left(\frac{x}{\ell}\right)$, for $x$  prime to  $\ell$,  where $\ell$  is an odd prime (where $(\frac{.}{.})$ is the Legendre-Jacobi symbol). If $x$ is odd, let us define $\chi_{-4}(x) = \left(\frac{-4}{x}\right) = (-1)^{(x-1)/2}$ and $\chi_8(x)=\left(\frac{8}{x}\right)=(-1)^{(x^2-1)/8}$. Given an integer $\delta,$ let 
$$
X(\delta) \;=\; \{\chi_{\ell}:\ell\textrm{ odd prime with }\ell\mid\delta\} \;\cup\; \{ \chi_{-4},\chi_{8}, \chi_{-4}\chi_{8}\}
$$
be the indicated set of  characters. If $f$ is a primitive ternary form (as in Brandt \cite{brandt1951zahlentheorie}) or  $f$ is an improperly primitive ternary form (as in Smith \cite{smith}) and if $\chi_{n}\in X(\delta),$ for some $\delta,$ then $\chi_n$ is called an \textit{assigned character} of $f$  if the following relation holds:  
\begin{equation}
\label{assignedcharactersdefinition}
    \chi_{n}(r_1) \;=\; \chi_{n}(r_2), \text{ for any } r_1,r_2\text{ represented by }f\text{ with }\gcd(r_i,n)=1.
\end{equation} 
This common value is denoted by $\chi_{n}(f).$ In particular, when this holds for some $\chi_n\in\{ \chi_{-4},\chi_{8}, \chi_{-4}\chi_{8}\},$  we say that $\chi_n$ is a \textit{supplementary assigned character} of the form $f.$ When both $\chi_{-4}$ and $\chi_{8}$ are supplementary assigned characters of the form $f,$   it is easy to see that Eq.~(\ref{assignedcharactersdefinition}) also holds for  $\chi_{-4}\chi_8.$ Since it is customary to list only two of these characters in this case, we usually do not list $\chi_{-4}\chi_8$ explicitly in this situation.

Let $X(f)$ denote the set of the assigned characters of the primitive ternary form $f$  with the basic invariants $I_1$ and $I_2$, and let $F:=F_f^B$. If we let $X^*(\delta)=X(\delta)\setminus\{\chi_8,\chi_{-4}\chi_{8}\}$, then we see by Brandt \cite{brandt1951zahlentheorie}, \S19, that
\begin{align}
\label{eq1: set_assigned_char_brandtI1}
    &X(f) \;=\; X(I_1) \   \text{ if }\  32 \mid I_1, \text{ and }  X(f) \;=\; X^*(I_1) \ \text{ if } \ I_1 \equiv 16\Mod{32},\\
    \label{eq2: set_assigned_char_brandtI2}
     &X(F) \;=\; X(I_2)\  \text{ if } \  32\mid I_2, \text{ and }  X(F) \;=\; X^*(I_2)\  \text{ if } \  I_2 \equiv 16\Mod{32}.
\end{align}
Recall also from \cite{brandt1951zahlentheorie}, loc.\ cit., that $f$ does not have a supplementary assigned character if $I_1(f)$ is odd.
From this discussion and Eq.~(\ref{I1Omega}), we can see that for  an improperly primitive positive definite ternary form $f$, the set of supplementary assigned characters $X_s(F_f)$ of $F:=F_f$ is
\begin{equation}
\label{eq: assigned_character_F_f_Smith_sense}
    X_s(F) \;=\; \{\chi_{-4}, \chi_8\} \  \text{ if } \  4\mid \Delta_f, \text{ and }     X_s(F)=\{\chi_{-4}\}  \ \text{ if } \ \Delta_f\equiv 2 \Mod{4}
\end{equation}
 as was explicitly shown in Table I, Case B of \cite[p.~459]{smith}. Also, by the same table, we have that
\begin{equation}
\label{Smithcharacter}
    \chi_{-4}(F_f) \;=\; -\chi_{-4}(\Omega_f).
\end{equation} 

\section{The refined Humbert invariant} \label{refinedhumbert}
In this section, we recall the concept of the refined Humbert invariant which was introduced by Kani \cite{kani1994elliptic}, \cite{MJ}. This is a quadratic form that is intrinsically attached to a principally polarized abelian surface. It is a useful ingredient since it can be used to translate  geometric problems into arithmetic ones, and thus it can be used to solve many interesting geometric problems; cf.\ \cite{MJ}, \cite{kani2014jacobians}. 

Let $A$ be an abelian surface over an algebraically closed field $K$ with $\charec K=0,$ and assume that $A$ has a principal polarization $\lambda: A\DistTo\hat{A}.$ Thus, $\lambda=\phi_{\theta}$ for a (unique) $\theta\in \NS(A)=\Div(A)/\equiv,$ where $\equiv$ denotes numerical equivalence; (see \cite[p.~60]{mumford1970abelian} for the discussion of $\phi$).  Here  $\theta$ is the (numerical) equivalence class $\cl(\Theta)$ of an ample divisor $\Theta$ in $\Div(A).$ Since we have by the Riemann-Roch Theorem (cf.\ \cite[p.~150]{mumford1970abelian}) that $4\deg(\phi_{\theta})=(\theta.\theta)^2$ for any ample $\theta,$ it follows that $(\Theta.\Theta)=2.$ For a principally polarized abelian surface $(A,\theta),$ let us put 
$$
\Tilde{q}_{(A,\theta)}(D) \;=\; (D.\theta)^2-2(D.D), \textrm{ for  } D\in \NS(A),
$$ where $(.)$ denotes the intersection number of divisors. By  \cite[p.~200]{kani1994elliptic}, the form $\Tilde{q}_{(A,\theta)}$ induces a positive definite quadratic form $q_{(A,\theta)}$ on the quotient module 
$$
\NS(A,\theta):=\NS(A)/\mathbb{Z}\theta.
$$
The quadratic form $q_{(A,\theta)}$ or, more correctly, the quadratic module  $(\NS(A,\theta), q_{(A,\theta)})$ is called the \emph{refined Humbert invariant} of the principally polarized abelian surface $(A,\theta)$; see \cite{MJ}, \cite{kani2014jacobians}. Let $\rho=\rank(\NS(A))$ be the Picard number of $A$. Thus $q_{(A,\theta)}$ gives rise  to an equivalence class of integral quadratic forms in $\rho-1$ variables since $\NS(A,\theta)\simeq\mathbb{Z}^{\rho-1}.$ Hence, when $\rho=4,$ the (corresponding) refined Humbert invariant is actually a ternary  quadratic form. 

\begin{remark} \label{rem: ref_hum_inv_relations_literature}
   It is interesting to observe that many authors have worked on geometric structures of the elements $(A,\theta)$ of $\mathcal{A}_2(K)$ by introducing a quadratic form that is equivalent to the refined Humbert invariant without mentioning or recognizing this equivalence in general.  For instance,  Hashimoto \cite{Hashimoto},  Runge \cite{Runge},  Lin and Yang \cite{quaternionic} conducted research in this manner for the particular case where the endomorphism algebra of the abelian surface $A$ is an indefinite quaternion algebra over $\Q$. The equivalence between the quadratic forms introduced or studied by these authors and the refined Humbert invariant will be shown in the author's PhD thesis.   Especially, Hashimoto, Lin and Yang obtained their associated quadratic form to $(A,\theta)$ by considering the set of `singular relations' in Humbert's sense. In fact, their approach has already been studied in \cite{kani1994elliptic} in the most general settings, where the refined Humbert invariant is introduced. By using the results in \S5 of \cite{kani1994elliptic},  we can observe that the associated quadratic form considered by these authors is equivalent to the refined Humbert invariant. This means that one can use Theorem \ref{maintheorem} and Theorem 1 of \cite{refhum} to characterize \textit{all} CM points on the intersection of two (different) Shimura curves considered in the papers of those authors.

  Additionally, R. Auffarth \cite{auffarth} gave a generalization
of the refined Humbert invariant for principally polarized abelian varieties to characterize abelian
varieties that contain an elliptic curve (by stating explicitly its relation with the refined Humbert invariant).
\end{remark}



Let $\mathcal{P}(A)\subset\NS(A)$ denote the set of \textit{principal polarizations} of $A.$ Thus, we see as in \cite[P.~141]{kani2014jacobians} that
$$
\mathcal{P}(A) \;=\; \{\cl(D)\in\NS(A): D\in \Div(A) \text{ is ample and } (D.D)=2\}.
$$
Since we aim to classify the imprimitive ternary forms which are equivalent to some refined Humbert invariant, we actually work with a subset of $\mathcal{P}(A).$ Before discussing it,  
 for two isogenous elliptic curves $E_1$ and $E_2$, let $q_{E_1,E_2}$ denote the \textit{degree map}  on $\Hom(E_1,E_2)$ which is defined by $q_{E_1,E_2}(h)=\deg(h),$ for $h\in \Hom(E_1,E_2).$

\begin{proposition}
\label{imprimitiveform}
 If $q_{(A,\theta)}$ is a ternary form, then $A\simeq E_1\times E_2$ is a product abelian surface, where $E_1$ and $E_2$ are isogenous CM elliptic curves. Moreover, we have that \begin{equation}
 \label{determinantoftherefined}
     \det(q_{(A,\theta)}) \;=\; 32\det(q_{E_1,E_2}), \text{ and so, } \disc(q_{(A,\theta)}) \;=\;  16\disc(q_{E_1,E_2})
 \end{equation}
 \end{proposition}
 \begin{proof}
 The  first statement follows as in the proof of Theorem 1 of \cite{refhum}. Indeed, since $q_{(A,\theta)}$ is a ternary form, it follows that $\rank(\NS(A))=4$ (as was discussed above). Thus, by the structure theorems for $\End(A)$ (cf.\ Proposition IX.1.2 of \cite{van2012hilbert}), it follows that  $A\sim E\times E,$  for some CM elliptic curve $E.$ By Shioda and Mitani's Theorem \cite{MitaniShioda} (or by Theorem 2 of \cite{kani2011products}), the first statement follows. The first equation of (\ref{determinantoftherefined}) follows from Lemma 30 of \cite{ESCI} together with Eq.~(37) of \cite{ESCI}. The second equation of (\ref{determinantoftherefined}) follows from the first equation of (\ref{determinantoftherefined}) and \cite[p.~2]{watson1960integral}.
  \end{proof}

We present a useful numerical classification of the set $\mathcal{P}(A)$ when $A$ is a product surface. To do so, we first recall some useful results and notations from \cite{MJ}. In particular, we recall the type of the binary forms studied in \cite{MJ}. These forms have important connections with the ternary forms that are the subject of this article.

 To state the result related to these binary forms, given an integer $\delta\geq1,$ let
 $$
 P(\delta) \;:=\; \{(n_1,n_2,k)\in\Z^3 : n_1,n_2>0  \textrm{ and } n_1n_2-k^2\delta=1 \},
 $$
 and for $s=(n_1,n_2,k)\in P(\delta),$ let 
 \begin{equation}
 \label{q_s}
    q_s(x,y) \;:\; =n_1^2x^2+2k\delta(n_1n_2+2)xy+n_2^2\delta(n_1n_2+3)y^2.
\end{equation} By \cite[p.~27]{MJ}, its discriminant  is \begin{equation}
\label{discq_s}
    \disc(q_s) \;=\; -16\delta.
\end{equation}   

We state one more result related to $q_s$ when it is an imprimitive binary form.  For this, put
$$
P(\delta)^{\ev} \;:=\; \{(n_1,n_2,k)\in P(\delta) : 2\mid\gcd(n_1,n_2) \}.
$$  By Theorem 13 of \cite{MJ} (with its proof), we obtain the following characterization of the imprimitive binary forms $q_s.$
\begin{proposition} \label{classificationoftypedelta}
Let $\delta\geq 1,$ and let $q$ be a binary form. Then $q=4q',$ for some primitive form $q'$ of discriminant $-\delta\equiv 1\Mod{4}$ which is in the principal genus $\gen(1_{-\delta})$ if and only if $q\sim q_s,$ for some $s\in P(\delta)^{\ev}.$
\end{proposition}

If $A=E_1\times E_2$ is a product surface, where $E_1$ and $E_2$ are two elliptic curves, then we have an isomorphism by Proposition 23 of \cite{MJ}:
\begin{equation}
\label{eq: isomorphism_NS(A)}
    \mathbf{D} \;=\; \mathbf{D}_{E_1,E_2}:\Z \oplus\Z\oplus\Hom(E_1,E_2)\DistTo \NS(A).
\end{equation}

For an abelian surface $A$, let  $\mathcal{P}(A)^{\odd}:=\mathcal{P}(A)\setminus\mathcal{P}(A)^{\ev}$, where 
$$
\mathcal{P}(A)^{\ev} \;=\; \{\theta\in\mathcal{P}(A):2\mid  (D.\theta),\forall D\in \NS(A)\}.
$$
 By  \cite{kani2014jacobians}, we have the following useful result for the refined Humbert invariants.
 
 \begin{proposition}
 \label{Q2aoddcase}
Let $A\simeq E_1\times E_2$ be a product surface with CM elliptic curves $E_1\sim E_2$, and let $f_{q}:=x^2+4q_{E_1,E_2}$. Then 
\begin{enumerate}
\item[\emph{(i)}] if $\theta\in \mathcal{P}(A),$ and $p$ is an odd prime or $p=\infty,$ then $q_{(A,\theta)}\sim_pf_q,$
\item[\emph{(ii)}] if $\theta\in \mathcal{P}(A)^{\odd}$, then  $q_{(A,\theta)}\in\gen(f_q).$ 
\end{enumerate}
\end{proposition}

\begin{proof}
The first assertion (i) follows from Corollary 19 of \cite{kani2014jacobians} because $f_q=q_{(A,\theta)},$ where $\theta=\textbf{D}(1,1,0)$. Thus, if we apply Theorem 20 of \cite{kani2014jacobians} to the quadratic module $(\NS(A), q_A),$ where $q_A$ is the intersection pairing as defined in \cite[p.~142]{kani2014jacobians}, then we get that  $q_{(A,\theta)}\in\gen(f_q),$ and so the second assertion (ii) follows.
\end{proof}

This result leads to significant formulas that will be discussed below; cf.\ Corollary \ref{cor: characters_for_ref_hum}. This result also gives a useful relation what we state now.
 \begin{corollary}
 \label{imprimitiveiffeven}
 If $\theta\in\mathcal{P}(A),$ where $A$ is as in \emph{Proposition \ref{Q2aoddcase}}, then 
 \begin{equation*}
     q_{(A,\theta)} \text{ is an imprimitive ternary form }  \ \Leftrightarrow \ \theta\in\mathcal{P}(A)^{\ev}.
 \end{equation*}
 \end{corollary}
 
\begin{proof}
If we take $\theta\in\mathcal{P}(A)^{\odd},$ then we see that $q_{(A,\theta)}\in\gen(f_q),$ for some binary quadratic form $q$ by Proposition \ref{Q2aoddcase}(ii), and thus $q_{(A,\theta)}$ is primitive since clearly $f_q$ is primitive.

Conversely, let us take $\theta\in\mathcal{P}(A)^{\ev}$. We claim that 
\begin{align*}
    \label{eq: 4_divide_ref_hum}
    4 \;\mid\; q_{(A,\theta)}(D), \ \forall \  D \in \NS(A).
\end{align*}
To prove this, recall that $2 \mid (D.\theta)$, for any $D\in \NS(A)$ by definition since $\theta \in \mathcal{P}(A)^{\ev}$. Also, the self-intersection number of a divisor on $A$ is even by the Riemann-Roch Theorem (see e.g., \cite[p.~150]{mumford1970abelian}), and thus, the claim follows. From this, the assertion follows because it shows that $4\mid \cont(q_{(A,\theta)})$.
\end{proof}

Note that if $\theta\in\mathcal{P}(A),$ where $A$ is as in Proposition \ref{Q2aoddcase}, then we obtain similarly that $q_{(A,\theta)}$ is a primitive ternary form if and only if $\theta\in\mathcal{P}(A)^{\odd}.$ In this respect, the classification of the ternary quadratic forms which are equivalent to some refined Humbert invariant $q_{(A,\theta)}$ for some $\theta\in\mathcal{P}(A)^{\odd}$ was done by Kani \cite{refhum} as was mentioned in the introduction. This is the \textit{primitive ternary case.} The main purpose of the present article is to give a similar classification for the ternary quadratic forms which are equivalent to some refined Humbert invariant  $q_{(A,\theta)}$ for $\theta\in\mathcal{P}(A)^{\ev}$, so this is the \textit{imprimitive ternary case} (cf.\ Corollary \ref{imprimitiveiffeven}).

 Let $A=E_1\times E_2$ be a product surface. While $\mathcal{P}(A)^{\odd}$ is always nonempty, $\mathcal{P}(A)^{\ev}$ may be an empty set for certain cases. 
 To discuss this, we use the following notations. If $q$ is a positive definite integral quadratic form, then let 
 $$
 R(q) \;=\; \{q(x)>0:x\in \Z^n\}
 $$
 denote the set of positive numbers represented by $q,$ and let 
 $$
 R_{a,m}(q) \;=\; \{r\in R(q): r\equiv a\Mod{m}\},
 $$
 for given integers $a\textrm{ and }m.$

\begin{proposition}
\label{eventhetaclassification}
Let $A=E_1\times E_2$ be a product surface with CM elliptic curves $E_1\sim E_2$, and let $d$ denote the discriminant of the form $q_{E_1,E_2},$ let $t$ denote the content of $q_{E_1,E_2}.$ Then we have that
\begin{enumerate}
\item[\emph{(i)}] $\mathcal{P}(A)^{\ev}\neq\varnothing$ if and only if $4\mid d$ and $R_{3,4}(q_{E_1,E_2})\neq\varnothing.$
\item[\emph{(ii)}] If $d$ is odd or if $t$ is even, then  $\mathcal{P}(A)^{\ev}=\varnothing$.
\end{enumerate}
\end{proposition}

\begin{proof}
We apply Proposition 32 of \cite{kani2014jacobians} to the quadratic module
$(\NS(A), q_A).$ Then Eq.~(6) of \cite{kani2014jacobians} shows that the hypothesis of that proposition
holds and that $q'=q_{E_1,E_2},$ where $q'$ is as defined in that proposition.
Then the first assertion (i) follows from Proposition 32 and Remark 33 of \cite{kani2014jacobians}. Also, (ii) easily follows  from (i).
\end{proof}


For later use, we prove the following useful result here that describes a `canonical' principal polarization $\theta \in \mathcal{P}(A)^{\ev}$ (if this is non-empty).
\begin{proposition} 
\label{eventhetaprop}
Let $A=E_1\times E_2$ be as in \emph{Proposition \ref{eventhetaclassification}}. If $\theta= \emph{\textbf{D}}(a,b,h)$ is an element of $\mathcal{P}(A),$  then we have that 
\begin{equation}
\label{eventheta}
    \theta \;=\; \emph{\textbf{D}}(a,b,h)\in\mathcal{P}(A)^{\ev} \ \Leftrightarrow\ 2\mid a,\ 2\mid b,\ 4\mid\disc(q_{E_1,E_2}).
    \end{equation}
    If this is the case, then $\deg(h) \in R_{3,4}(q_{E_1,E_2})$. 
    
    In addition,  if $4\mid \disc(q_{E_1,E_2})$ and if there is an $h\in\Hom(E_1,E_2)$ such that $r:=q_{E_1,E_2}(h)\equiv3\Mod{4},$ then 
    $$
    \theta':= \emph{\textbf{D}}(2, (r+1)/2,h)\in\mathcal{P}(A)^{\ev}.
    $$
\end{proposition}
 \begin{proof}
 As in the proof of Proposition \ref{eventhetaclassification}, if we apply Proposition 32 of \cite{kani2014jacobians} to the quadratic module $(\NS(A), q_A),$ then we can take $q'=q_{E_1,E_2},$ and so we see 
  from Eq.~(30) of \cite{kani2014jacobians} that $$
  \textbf{D}(a,b,h)\in\mathcal{P}(A)^{\ev} \;\Leftrightarrow\; 2\mid a, \ 2\mid b,\textrm{ and } \beta_{q'}(h,h')\equiv 0\Mod{2}, \ \forall h'\in\Hom(E_1,E_2),
  $$
   where $\beta_{q'}$ is the bilinear map associated to $q';$ cf.\ Eq.~(13) of \cite{kani2014jacobians}. By Remark 33 of \cite{kani2014jacobians}, we see that $ 2\mid a,\ 2\mid b,\ 4\mid \disc(q_{E_1,E_2}) \Leftrightarrow 2\mid a,\ 2\mid b,\ \beta_{q'}(h,h')\equiv 0\Mod{2},$ for all $h'\in\Hom(E_1,E_2),$ and thus Eq.~(\ref{eventheta}) holds.
  
If $\theta = \textbf{D}(a,b,h)\in \mathcal{P}(A)^{\ev}$, then $ab\equiv 0\Mod{4}$ by (\ref{eventheta}). Then we obtain by  Corollary 25 of \cite{MJ} that $-q_{E_1,E_2}(h)\equiv 1\Mod{4}$, and so $r=q_{E_1,E_2}(h)\in R_{3,4}(q_{E_1, E_2})$, as desired. Also, note that $\theta'\in\mathcal{P}(A).$ Indeed, since $2(({r+1})/{2})-r=1,$ we have that $\theta'\in\mathcal{P}(A)$ by Corollary 25 of \cite{MJ} again. So, this part easily follows from (\ref{eventheta}). 
\end{proof}

\section{Quadratic forms in \texorpdfstring{$\Theta(A)^{\ev}$}{}}
\label{section: quadratic forms in Q(A)}
 

Observe that  Proposition \ref{Q2aoddcase}(ii) gives that $\Theta(A)^{\odd} := \{q_{(A,\theta')}: \theta'\in\mathcal{P}(A)^{\odd}\}/\! {\sim}$ lies in a \textit{single} genus. In this section, we study the set 
$$
\Theta(A)^{\ev} \;:=\; \{q_{(A,\theta')} : \theta'\in\mathcal{P}(A)^{\ev}\}/\! {\sim},
$$
and we see in Theorem \ref{cor: Theta_lies_genus} below that the set  $\Theta(A)^{\ev}$ also lies in a single genus in most cases, but not in all cases. This theorem is the main result of this section, and it will be a key tool for the main results of this article. More precisely, the sufficient part of our classification theorem immediately follows from this theorem and Theorem \ref{theorem1}.

The study of the representation of a quadratic form by another quadratic form is a key ingredient in this section. If a binary integral quadratic form $\phi$ is \textit{primitively represented} by a ternary integral quadratic form $f$ (as defined in \cite[p.~25]{dicksonsbook}), we write $f\rightarrow\phi.$ If a quadratic form $f$ \textit{primitively represents} an integer $n$ (as defined in \cite[p.~8]{watson1960integral}), we  write $f\rightarrow n.$
Throughout the article, we use many useful results in \cite{MJ}. In that paper, a slightly different definition of the symbol $f\rightarrow\phi$ was used. To recall this definition, 
let us consider 
two quadratic $\mathbb{Z}$-modules $(\Z^3,f)$ and $(\mathbb{Z}^2,\phi).$ If there exists an injective homomorphism $\psi:\mathbb{Z}^2\xhookrightarrow{} \Z^3$ with $\Z^3/\psi(\mathbb{Z}^2)$ torsion-free such that $\phi=f\circ \psi,$ then we say that $f$ \textit{primitively represents} the quadratic form $\phi$. Note that the invariant factor theorem shows that both definitions are equivalent. 

\subsection{The properties of the forms in \texorpdfstring{$\Theta(A)^{\ev}$}{}}
In this subsection, we provide the properties of the forms in $\Theta(A)^{\ev}$. Recall that we have established the relations between the basic invariants $I_1, I_2$ and $\Omega, \Delta$ and between the reciprocals $F$ and $F^B$ in Proposition \ref{F=F^B}.  Now, we explicitly determine them for the forms in $\Theta(A)^{\ev}$. The following result from E. Kani's handwritten notes was shared with the author. 

\begin{proposition}
\label{invariants}
Let $A=E_1\times E_2$ be a product surface with CM elliptic curves $E_1\sim E_2$, and let $\theta \in \mathcal{P}(A)^{\ev}$. Let us put
     $f:=\frac{1}{2}q_{(A,\theta)},$ $d:=-\disc(q_{E_1,E_2}), \ t:=\cont(q_{E_1,E_2}) \textrm{ and } d':=d/t^2.$
 Then $f$ is an improperly primitive form and basic invariants of $f/2$ are   
 \begin{equation}
    I_1(f/2)\;=\;-t \textrm{ and } I_2(f/2)\;=\;-4d', \textrm{ and thus, } \disc(f/2)\;=\;\frac{-t^2d'}{4}. \label{I1I2basicinvariants} 
\end{equation} Therefore, the following equations hold 
\begin{align}
    \label{basicinvariantsDicksonsense}
     \Omega_f &\;=\;t \ \textrm{ and } \ \Delta_f\;=\;d'/2,\ \text{ and } \\
      \label{todd4dividesd'}
      t & \ \text{ is odd and } \ 4\mid d'.
\end{align}
Moreover, if $\theta = \emph{\textbf{D}}(n_1,n_2,kh)\in\mathcal{P}(A)^{\ev}$, for a primitive $h\in\Hom(E_1,E_2)$, then 
  $q_{E_1,E_2}(h)/t$ is represented by the reciprocal $F_f=F^B_{f/2}$ of $f/2.$
\end{proposition}

\begin{proof}
Firstly, observe that $\cont(q_{(A,\theta)})=4.$ Indeed, it is easy to see that $4\mid q_{(A,\theta)}(D),$ $\forall D\in \NS(A)$ (as in the proof of Corollary \ref{imprimitiveiffeven}) and so, it follows that $4\mid\cont(q_{(A,\theta)}).$  Since $\theta\in\mathcal{P}(A)^{\ev},$ we have that $n_1$ and $n_2$ are even by condition (\ref{eventheta}), and $n_1,n_2>0$ and $n_1n_2-\deg(kh)=1$ by Corollary 25 of \cite{MJ}. Hence, this gives that $s=(n_1,n_2,k)\in P(\deg(h))^{\ev}.$ By Proposition 29 of \cite{MJ}, we have that $q_{(A,\theta)}\rightarrow q_s,$ where $q_s$ is as in equation (\ref{q_s}). Since $s\in P(\deg(h))^{\ev},$ we have that $\cont(q_s)=4$ by Proposition \ref{classificationoftypedelta}. Thus, $\cont(q_{(A,\theta)})\mid 4,$ and so, we have that $\cont(q_{(A,\theta)})=4,$ which proves the observation.

Since $\theta\in\mathcal{P}(A)^{\ev}$, it follows that $\deg(h)$ is odd. Since $s\in P(\deg(h))^{\ev}$, we see by equation (\ref{discq_s}) that $\disc(q_s)=-16\deg(h),$ so $\disc(\frac{1}{4}q_s)=-\deg(h)$ is odd. Hence, it follows that $\frac{1}{2}q_s$ is improperly primitive. Since $f=\frac{1}{2}q_{(A,\theta)}\rightarrow \frac{1}{2}q_s,$ it follows that $f$ is also improperly primitive. Indeed, we know that $f/2$ is a primitive form. If $f$ were not improperly primitive, then we would have that all the non-diagonal coefficients of $f$ would be divisible by 4. If we consider the coefficient matrices $A(f)$ and $A(\frac{1}{2}q_s)$ of the forms $f$ and $\frac{1}{2}q_s$, (see \cite[p.~2]{watson1960integral} for the definition of the coefficient matrix), then there exists an integral $3\times 2$ matrix $T$ such that $T^tA(f)T=A(\frac{1}{2}q_s)$ since $f\rightarrow \frac{1}{2}q_s.$ Since all the non-diagonal coefficients of $f$ are divisible by $4$ and $\cont(f)=2,$ it follows that all entries of $A(f)$ are divisible by $4,$ and so all entries of the $A(\frac{1}{2}q_s)$ are divisible by $4.$ But this would imply that the non-diagonal coefficient of the $\frac{1}{2}q_s$ would be divisible by $4$, which is not possible since it is improperly primitive.  This proves that $f$ is improperly primitive.

Since $f$ is improperly primitive, $I_1(f/2)$ is a negative odd number by Eq.~(\ref{I1odd}). We know that $2f\sim_p f_q,$ for all odd primes $p,$ where $f_q=x^2+4q_{E_1,E_2},$ by Proposition \ref{Q2aoddcase}, and thus, $f/2\sim_p f_q,$ for all odd primes $p$ since $2$ is unit in $\Z_p$. Since $f/2=\frac{1}{4}q_{(A,\theta)}$ and $f_q$ are primitive forms and $p$-adically equivalent for all odd primes $p$, we get that  $v_p(I_1(f/2))=v_p(I_1(f_q))$ for odd primes $p$, where $v_p(n)$ denotes the exponent of the largest prime power of $p$ that divides $n$. Since $I_1(f_q) = -16t$ by the proof of Proposition 53 of \cite{kani2014jacobians}, and $t$ is odd (cf.\ Proposition \ref{eventhetaclassification}(ii)), we thus see that $I_1(f/2)=-t.$ 

 By Eq.~(\ref{determinantoftherefined}), we see that 
 $$
 \disc(f/2)\;=\; \disc\left(\frac{1}{4}q_{(A,\theta)}\right)\;=\; \frac{1}{4^3}\disc(q_{(A,\theta)})\;=\;\frac{-d}{4}.
 $$
 Thus, we see that $I_2(f/2)=16(-d/4)/t^2=-4d'$ by Eq.~(\ref{discf}), so Eq.~(\ref{I1I2basicinvariants}) holds. Then Eq.~(\ref{basicinvariantsDicksonsense}) just follows from Eq.~(\ref{I1I2basicinvariants}) together with Eq.~(\ref{I1Omega}). Additionally, Eq.~(\ref{todd4dividesd'}) follows from Eq.~(\ref{I1odd}).

Lastly, since  $f\rightarrow \frac{1}{2}q_s$ and $\disc(\frac{1}{2}q_s)=-4\deg(h)$ by what was mentioned above, we have that $\det^{\text{D}}(\frac{1}{2}q_s)=-\deg(h)$ by Eq.~(\ref{dicksondeterminant}). So, it follows from Theorem 27 of \cite{dicksonsbook} that $\deg(h)/t\in R(F_f).$ Thus,  the last assertion follows because $F_f=F_{f/2}^B$ by Proposition \ref{F=F^B}. 
 \end{proof}

This result also shows that condition (\ref{classificationconditions1}) of a positive definite ternary form $f$  is necessary for $f$ to be equivalent to some (imprimitive) refined Humbert invariant $q_{(A,\theta)}$. 

A nice application of the classification of the binary forms $q_s$ is that it can be used to prove the necessary condition (\ref{classificationconditions2}) of the second main result.

\begin{proposition} 
\label{necessaryconditionforii}
If an imprimitive ternary form $f$ is equivalent to some refined Humbert invariant $q_{(A,\theta)},$ then there exists an integer $n$ relatively prime to $\disc(f)$ such that $f$ represents $(2n)^2$.
\end{proposition}

\begin{proof}
 By Proposition \ref{imprimitiveform}, we have that $A=E_1\times E_2$ is a product surface, where $E_1 \sim E_2$ are CM elliptic curves since $q_{(A,\theta)}$ is a ternary form. Let us put $\theta=:\textbf{D}(n_1,n_2,kh),$ where $h\in\Hom(E_1,E_2)$ is primitive. Since $\theta\in\mathcal{P}(A)^{\ev}$ by Corollary \ref{imprimitiveiffeven}, we have that $q_{(A,\theta)}\rightarrow q_s,$ where $s\in P(\deg(h))^{\ev}$ as in the proof of Proposition \ref{invariants}.


Since $s\in P(\deg(h))^{\ev},$ we have by Proposition \ref{classificationoftypedelta} that $q_s=4q,$ for some primitive form $q$ lying in the principal genus $\gen(1_{-\deg(h)})$. Also, we see that there exists an integer $n\geq 1$ with $\gcd(n,\disc(f))=1$ such that $ n^2R(1_{-\deg(h)})\subset R(q)$ by Proposition 14 of \cite{refhum} (since $q\in\gen(1_{-\deg(h)}))$. Since a principal form represents 1, it follows that $q$ represents $n^2$. Therefore, since $q_{(A,\theta)}\rightarrow q_s=4q$, we obtain that $4n^2\in R(q_{(A,\theta)})=R(f),$ which proves the assertion. 
\end{proof}

The following corollary is very useful, because it says that the values of all, with the possible exception of one, of the assigned characters of the imprimitive ternary forms $q_{(A,\theta)}$ and $q_{(A,\theta')}$ are equal. For a primitive ternary form 
$f$, let $X^*(f)$ denote the set $X(f)\setminus\{\chi_8, \chi_{-4}\chi_8\}$. 

\begin{corollary}
    \label{cor: characters_for_ref_hum}
    Let $A=E_1\times E_2$ be as in \emph{Proposition \ref{invariants}}. Let $\theta_i \in\mathcal{P}(A)^{\ev}$,  and let $4f_i:=q_{(A,\theta_i)}$, for $i=1,2$. Then $f_1$ and $f_2$ have the same basic invariants $I_1$ and $I_2$, and
$4f_1 \sim_p 4f_2$, for all odd primes $p$. Moreover,
\begin{equation}
\label{cor: eq_characters}
     1 \;=\; \chi(f_1) \;=\; \chi(f_2),\ \forall\chi\in X(f_1), \textrm{ and }\  \chi(F_{f_1}^B)\;=\;\chi(F_{f_2}^B),\ \forall\chi\in X^*(F_{f_1}^B).
\end{equation}
\end{corollary}

\begin{proof}
We have that $\cont(q_{(A,\theta_i)}) =4$ by Proposition \ref{invariants}, for $i=1, 2,$ and also we see from this proposition that $I_1 := I_1(f_1) = -\cont(q_{E_1,E_2}) = I_1(f_2)$ and  $I_2:= I_2(f_1) = -4\disc(q_{E_1,E_2})/I_1^2= I_2(f_2)$, and this proves the first assertion. 

Since $\cont(q_{(A,\theta_i)}) =4$, we see that $q_{(A,\theta_i)}\sim_{p} \frac{1}{4} q_{(A,\theta_i)}$ for all odd primes $p$ (as in the proof of Proposition \ref{invariants}), and so $f_1\sim_p f_2$, for all odd primes $p$ by Proposition \ref{Q2aoddcase}(i), and so $4f_1 \sim_p 4f_2$, as stated.

We now show the first equation of (\ref{cor: eq_characters}).    Since the $2f_i$ are improperly primitive ternary forms by Proposition \ref{invariants},  $I_1$ is odd by (\ref{I1odd}). Hence, we have that $X(f_2) = X(f_1)=\{\chi_{p}: p\mid I_1,\ p \text{ prime}\},$ and that $f_i$ has no  supplementary assigned characters (cf.\ \S\ref{ternaryforms}). Therefore, since $f_1 \sim_p f_2 \sim_p f_q$, for any prime $p\mid I_1$ (where $f_q(x,y,z)=x^2+4q_{E_1,E_2}(y,z)$ as before)  and since $I_1(f_q) = 16I_1$ (cf.\ Proposition \ref{Q2aoddcase} and the proof of Proposition \ref{invariants}), we obtain that
$$
\chi(f_1) \;=\; \chi(f_2) = \chi(f_q) \;=\; 1, \quad \forall \chi \;\in\; X(f_1).
$$
The last equality follows from the fact that $f_q(1,0,0)= 1$, and so this proves the first equality of (\ref{cor: eq_characters}).  

To show the second equality of (\ref{cor: eq_characters}), note that $X(F_{f_1}^B) = X(F_{f_2}^B)$ since the primitive forms $f_i$ have the same basic invariant $I_2$. Since  $16 \mid I_2=-8\Delta_{2f_i}$ by Proposition \ref{F=F^B}, it follows that $\chi_{-4} \in X(F_{f_1}^B)$; see (\ref{eq: assigned_character_F_f_Smith_sense}). We know that $F_{f_i}^B = F_{2f_i}$ by Proposition \ref{F=F^B}, and that $\chi_{-4}(F_{2f_i}) = - \chi_{-4}(\Omega_{2f_i}) $  by Eq.~(\ref{Smithcharacter}), for $i=1,2.$ Therefore, we obtain that $\chi_{-4}(F_{f_1}^B) = \chi_{-4}(F_{f_2}^B)$ since $\Omega_{2f_1} = \Omega_{2f_2}$ by Eq.~(\ref{I1Omega}).

 Since  $f_1 \sim_{p} f_2$ for all odd primes $p$, it follows  that $F_{f_1}^B \sim_{p} F_{f_2}^B$, and thus $\chi_p(F_{f_1}^B)=\chi_p(F_{f_2}^B)$, for all odd primes $p\mid I_2$, and so, this proves the second equation of (\ref{cor: eq_characters}), and hence the assertions follow.
\end{proof}

By Theorem \ref{cor: Theta_lies_genus} below, we see that the set $\Theta(A)^{\ev}$ may not lie in a single genus, unlike in  the set $\Theta(A)^{\odd}$. The reason why this is not the case is due to the following technical lemma. 

\begin{lemma}
\label{representstwonumbers} Assume that $q$ is a positive definite binary quadratic form with $R_{3,4}(q)\neq\varnothing.$ In addition, assume that $\disc(q)\equiv 16 \Mod {32}.$ Then the form $q$ primitively represents two numbers $r_1$ and $r_2$ such that $r_1\equiv 3 \Mod 8$ and $r_2\equiv 7 \Mod 8.$ 
\end{lemma}

\begin{proof}
Since $q$ primitively represents a number $r_1\equiv 3 \Mod 4$ by the hypothesis, we have that $q\sim r_1x^2+bxy+cy^2$ for some $b,c\in\Z,$ by Lemma 2.3 of \cite{davidcox}. 
Since $16 \mid \disc(q) = b^2 - 4r_1c$, it follows that $b$ is even, and so put $b=2b_1.$ Then this gives that $b_1^2-r_1c\equiv 4\Mod{8}.$   To prove the statement, we distinguish some cases.
\spm

\noindent{\bf Case 1.} Assume that $r_1\equiv 3 \Mod 8.$  
\sps

\noindent By the hypothesis, we have $b_1^2\equiv 3c+4\Mod{8}.$ We know that the set of the residue classes of the squares mod $8$ is $\{0,1,4\}.$ So there are three cases for $b_1^2$ in mod $8.$ We can easily see that $q(1,1)$, $q(1,1)$ and $q(0,1)$ respectively give a number $r_2 \equiv 7 \Mod{8}$ in these cases. Hence, $q$ primitively represents a number $r_2 \equiv 7 \Mod{8}$ in these three subcases of Case 1.
\spm

\noindent{\bf Case 2.} Assume that $r_1\equiv 7 \Mod 8.$
\sps

\noindent 
We have that $b_1^2\equiv 7c+4\Mod{8}.$ So there are three cases for $b_1^2$ in mod $8$, namely $0, 1, 4$. We can again see that $q(1,1)$, $q(0,1)$ and $q(1,1)$ respectively give a number $r_2 \equiv 3 \Mod{8}$, and so the assertion follows in this case.
\end{proof}

\subsection{Construction of \texorpdfstring{$(A,\theta)$}{}}

Let $f$ be a positive definite ternary form satisfying conditions (\ref{classificationconditions1}) and (\ref{classificationconditions2}). The basic purpose of this section is to explicitly construct a principally polarized abelian surface $(A,\theta)$ such that $q_{(A,\theta)}$ is genus equivalent to $f$. Since our construction is explicit, it is suitable to write an algorithm to find such $(A,\theta)$ for a given $f$. This was roughly outlined in \cite{refhum} when $f$ is primitive. This will be explicitly satisfied in the author's PhD thesis for any ternary form $f$. This construction plays the key role in the proof of Theorem \ref{cor: Theta_lies_genus} below. Moreover, it is a part of the main result
because the sufficient part of Theorem \ref{maintheorem} follows from Theorem \ref{theorem1} by this construction.

The following lemma from Remark 41 of \cite{kani2011products} has a crucial role in the following constructions of principally polarized abelian surfaces. Note that this lemma was explicitly proved in Lemma 26(a) of \cite{refhum}.

\begin{lemma}
\label{lemma26}
Let $q$ be a positive definite binary quadratic form. Then there exist two isogenous CM elliptic curves $E_1$ and $E_2$ such that $q_{E_1,E_2}\sim q.$
\end{lemma}

For  given positive odd $\Omega$ and even $\Delta$ integers, suppose that there is an integer $a$ such that
\begin{equation}
    \label{eq: construction_a}
   \gcd(a,\Delta) \;=\; 1, \quad \left(\frac{-2 \Delta}{a}\right) \;=\; 1 \quad \text { and } \quad a \Omega \;\equiv\; 3 \Mod{4}.
\end{equation}
Thus, there is an integer $b$ such that $-2 \Delta \equiv b^{2} \Mod{a}$. By replacing $b$ with $a-b$, if necessary, we may assume that $b$ is even. Hence, since $\Delta$ is even, it follows that $-2 \Delta \equiv b^{2} \Mod{4a}$, and so we put
$$
q_{\Delta}(x, y) \;:=\; a x^{2} \;+\; b x y \;+\; \frac{b^{2}+2 \Delta}{4 a} y^{2} .
$$
We see that $q_{\Delta}$ is a positive definite binary integral  quadratic form of discriminant $-2 \Delta$. Since $\gcd(a, 2 \Delta)=1$ by the hypothesis, it follows that $\gcd(a, b)=1$, and hence $q_{\Delta}$ is primitive. We now consider the form
\begin{equation}
\label{eq: construction_q_Omega_Delta}
    q_{\Delta, \Omega}(x, y) \;:=\; \Omega{q}_{\Delta} \;=\; \Omega a x^{2} \;+\; \Omega b x y \;+\; \Omega\frac{b^{2}+2 \Delta}{4 a} y^{2},
\end{equation}
which is a positive definite binary integral quadratic form of discriminant $-2 \Delta \Omega^{2}$ with content $\Omega$. By Lemma \ref{lemma26}, there exist two isogenous CM  elliptic curves $E_{1}, E_{2}$ such that $q_{E_1, E_2} \sim q_{\Delta, \Omega}$. We consider the abelian surface $A := E_{1} \times E_{2}$.
It is clear that $4\mid \disc(q_{\Delta, \Omega})$. Since $q_{\Delta,\Omega}(1,0) =  a \Omega \equiv 3 \Mod{4}$ by the hypothesis (see (\ref{eq: construction_a})), there is a primitive isogeny $h \in \Hom(E_{1}, E_{2})$ such that $q_{E_{1}, E_{2}}(h) = a \Omega$, and so  it follows from Proposition  \ref{eventhetaprop} that
\begin{equation}
\label{eq: construction_theta_for_Omega_Delta}
    \theta_a \;:=\; \textbf{D}(2, (\Omega a + 1)/2, h) \;\in\; \mathcal{P}(A)^{\ev}.
\end{equation}
 Therefore, we have constructed a principally polarized abelian surface $(A, \theta_a)$ such that 
\begin{equation}
\label{eq: notation_f_Delta_Omega}
    f^{\theta_a}_{q_{\Delta,\Omega}}\;:=\; q_{(A, \theta_a)} \ \text{ is an imprimitive ternary form,}
\end{equation}
  see Corollary \ref{imprimitiveiffeven}. 

Let us state what we proved by this construction here:
\begin{proposition}
\label{p: construction_A_th_1}
    For given positive odd $\Omega$ and even $\Delta$ integers, suppose that there is an integer $a$ satisfying the conditions in \emph{(\ref{eq: construction_a})}. Then, there is a product abelian surface $A = E_1 \times E_2$, where the degree map $q_{E_1,E_2}$ is a binary quadratic form of discriminant $-2\Omega^2\Delta$ with content $\Omega$, and also there is a principal polarization $\emph{\textbf{D}}(2, (\Omega a + 1)/2, h) \in \mathcal{P}(A)^{\ev},$ for a primitive isogeny $ h \in \Hom(E_1,E_2).$
\end{proposition}

To make more useful this construction we now observe that for given basic invariants $\Omega_{f/2}$ and $\Delta_{f/2}$, where a positive definite ternary form $f$ satisfying conditions (\ref{classificationconditions1}) and (\ref{classificationconditions2})   we prove the existence of the number $a$ satisfying the conditions in (\ref{eq: construction_a}).

\begin{proposition}
\label{prop: existence_number_a}
   Suppose that a positive definite ternary form $f$ satisfies conditions \emph{(\ref{classificationconditions1})} and \emph{(\ref{classificationconditions2})}. Then $\chi(f/4) = 1$, for all $\chi \in X(f/4)$. 
   
   Let $\Omega=\Omega_{f/2}$ and $\Delta=\Delta_{f / 2}$ be the basic invariants of $f/2$. Then there exists an integer $a$ satisfying the conditions in \emph{(\ref{eq: construction_a})}, for $\Omega$  and $\Delta$.  In addition, we have that $\chi_{p}(a)=\chi_p(F_{f/2})$, for all odd primes $p\mid \Delta$, and that $\chi_{8}(a) = \chi_{8}(F_{f/2})$ if $4\mid \Delta$.
\end{proposition}

\begin{proof} 
By  condition (\ref{classificationconditions1}), $f/2$ is improperly primitive, and so it follows by Proposition \ref{F=F^B} that $\Omega = -I_1(f/4)$ is odd and $\Delta$ is even. Thus, the set $X(f/4)$ of the assigned characters  of $f/4$ is $X(f/4) = \{\chi_p :  p \text{ odd prime, } p\mid \Omega\}$.  Recall that by assumption (\ref{classificationconditions2}), there is an integer $m$ with $m^2=\frac{1}{4}f(x,y,z)$, for some $x, y, z$, and $\gcd(m,\disc(f))=1.$ This immediately proves the first assertion that $\chi(f/4)=\chi(m^2) = 1$, for all $\chi \in X(f/4)$ (since $\gcd(m,\Omega) =1$; cf.\ (\ref{discf})).

Now, let us put $g=\gcd(x,y,z)$ and $n=m/g,$ so $n^2$ is primitively represented by $f/4,$ and $\gcd(n,\disc(f))=1.$ Clearly, it follows from  Eqs.~(\ref{discf}) and (\ref{I1Omega}) that $\gcd(n, \Omega\Delta) = 1$.

We claim that $F:=F_{f / 2}$ represents a binary form $\phi$ whose discriminant is $-8 \Delta{n}^{2}$. Indeed, since $f/2=F_{F}$ primitively represents $2n^2,$ we can apply Theorem 38 of \cite{dicksonsbook} to $F_{F},$ and thus conclude that there is a binary form $\phi$ with $F\rightarrow \phi$ and $\det^\text{D}(\phi)=-\Omega_{F}2n^2.$ Recall from Theorem 5 of \cite{dicksonsbook} that $\Omega_F=\Delta$.  Hence, we see that  $\disc(\phi)=-4\Omega_{F}2n^2=-8\Delta n^2;$ cf.\ Eq.~(\ref{dicksondeterminant}), and so the claim follows. 

We next observe that $\cont(\phi)$ is either 1 or 4. For this, we apply Proposition 10 of \cite{refhum} to $F$ with $C=2 n^{2}$. Note first that $F=F_{f / 4}^{B}$ is properly primitive by Proposition \ref{F=F^B}. Since $\gcd(C, 2 \Omega \Delta)=2$, it follows from  Proposition 10 of \cite{refhum} that $\cont(\phi)=1$ or 4.

Let us put $\phi=\kappa \phi^{\prime}$, for a primitive form $\phi^{\prime}$, i.e., $\kappa=1$ or 4. We now show that $\phi^{\prime}$ represents an integer $a$ relatively prime to $\Omega \Delta$ such that $a \Omega \equiv 3 \Mod{4}$.
\spm

\noindent{\bf Case 1.} Suppose that $\kappa=1$.
\sps

\noindent Since $\Delta$ is even, $\chi_{-4}$ is an assigned character of $F$ (cf.\ (\ref{eq: assigned_character_F_f_Smith_sense})), and $\chi_{-4}(F)=-\chi_{-4}(\Omega)$ by (\ref{Smithcharacter}). Moreover, we have that  $\frac{\disc(\phi')}{-4}=2\Delta n^2\equiv0\Mod{4},$ so $\chi_{-4}$ is an assigned character of $\phi'$ by \cite[p.~55]{davidcox}. Since $\phi'$ is primitive, it represents a number $a$ relatively prime to $\Omega\Delta$ (cf.\ Lemma 2.25 of \cite{davidcox}), and thus $a\in R(\phi')\subset R(F),$ and hence  $\chi_{-4}(a)=\chi_{-4}(\phi')=\chi_{-4}(F) = -\chi_{-4}(\Omega)$, which shows that $\chi_{-4}(a\Omega) = -1$, equivalently, $a\Omega \equiv 3 \Mod{4}$, as requested.
\spm

\noindent{\bf Case 2.} Suppose that $\kappa=4$.
\sps

\noindent We obtain that $2 n^{2} \equiv \Delta+4 \Mod{8}$ by Proposition 10 of \cite{refhum} since $\phi$ is not primitive. Then, since $n$ is odd, we have $2 n^{2} \equiv 2 \Mod{8}$, and so $\Delta \equiv 2 \Mod{4}$. Thus, $\disc(\phi')=-8 \Delta n^{2} / 4^{2}=-\Delta/2 n^{2}$ is odd. Hence, by the fact
that a primitive form with odd discriminant represents infinitely many primes $p \equiv 1 \Mod{4}$ and also $p^{\prime} \equiv 3 \Mod{4}$; cf.\ Corollary 4 of \cite{refhum}, there is a number $a \in R(\phi')$ relatively prime to $\Omega\Delta$ such that $a\Omega \equiv 3 \Mod{4}$ again.

Therefore, $\phi^{\prime}$ represents a number $a$ relatively prime to $\Delta \Omega$ such that $a \Omega \equiv 3 \Mod{4}$ in both cases. Thus, we have that $\left(\frac{-8 \Delta n^{2} / \kappa^{2}}{a}\right)=1$, for $\kappa=1, 4$  by Lemma 2.5 of \cite{davidcox}. Notice that $\kappa^{2} \mid 8 \Delta$ (since $\Delta$ is even), for $\kappa=1$ or 4 , and so it follows that $\left(\frac{-2 \Delta}{a}\right)=1$. Hence, there is a number $a$ satisfying the conditions in (\ref{eq: construction_a}), which proves the second assertion.

We know that $\chi_p$ is an assigned character of $F$, for any odd prime $p\mid \Delta$. Note that $\kappa^2 a \in R(\kappa\phi') = R(\phi) \subset R(F)$. Clearly, $\gcd(a\kappa, p) =1$, for any odd prime $p\mid \Delta$. Thus, $\chi_p(F) = \chi_p(\kappa^2 a) = \chi_p(a)$, for any odd prime $p\mid \Delta$. This proves the third assertion. Moreover, if $4\mid \Delta$, then $\kappa = 1$; otherwise, as was discussed in Case 2, $\Delta \equiv 2 \Mod{4}$. This implies that $a \in R(F)$  in this case, and so, since $\chi_8$ is an assigned character of $F$ (cf.\ (\ref{eq: assigned_character_F_f_Smith_sense})), $\chi_8(F) = \chi_8(a)$. Thus, the last assertion also follows.
\end{proof}

The following technical lemma is required for the constructions of the suitable principal polarizations; see (\ref{basicpp}) below. 

\begin{lemma}
\label{representsrelatively}
Let $q$ be a primitive binary form for which $\disc(q)=d$ is even. If $q$ primitively represents an odd number $a,$ then $q$ primitively represents a number $a'\equiv a\Mod 4$ with $(a',d)=1.$ In addition, if $16\mid d$, then $a'\equiv a\Mod 8$.
\end{lemma}

\begin{proof}
To prove the assertions we refine the argument of the proof of Proposition 4.2 of \cite{Duncanbinaryforms}. We may assume by Lemma 2.3 of \cite{davidcox} that $q=ax^2+bxy+cy^2$ with $b^2-4ac=-d$ since $q\rightarrow a$. Consider the following sets of prime divisors of $d:$ 
$$
\mathcal{P}_1 \;:=\; \{p\mid \gcd(a,c,d)\},\quad\mathcal{P}_2\;:=\;\{p\mid \gcd(a,d),p\nmid c\},
$$
$$
\mathcal{P}_3\;:=\;\{p\mid \gcd(c,d),p\nmid a\} \ \text{ and }\ \mathcal{P}_4\;:=\;\{p\mid d, p\nmid a, p\nmid c\}.
$$
It is clear that the union $\cup_i\mathcal{P}_i$ is the set of the prime divisors of $d,$ and that the $\mathcal{P}_i$'s are disjoint.

Put $x_i:=\prod_{p\in{\mathcal{P}_i}}p,$ and put $a':=q(x_2,x_3x_4),$ so $q\rightarrow a'$ since we have that $\gcd(x_2,x_3x_4)=1$ by the construction. We first show that $\gcd(a',d)=1.$ Assume not, so there is a prime number $p\mid\gcd(a',d).$  Hence, $p\in\mathcal{P}_i,$ for some $1\leq i\leq 4.$ It is also clear that $p\mid b$ for $p\in\mathcal{P}_i,$ with $i\in\{1,2,3\}$ since $b^2-4ac=d.$ To get a contradiction we need to distinguish four cases.
\spm

\noindent{\bf Case 1.} Assume that $p\in\mathcal{P}_1.$ Since $p\mid b,$ this is not possible since $q$ is primitive.
\spm 

\noindent{\bf Case 2.} Assume that $p\in\mathcal{P}_2.$ Write $a'\;=\;q(x_2,x_3x_4)=ax_2^2+bx_2x_3x_4+cx_3^2x_4^2$. Since $p\mid x_2$ and $p\mid a',$ it follows that $p\mid cx_3^2x_4^2.$ But, $p\nmid c,$ so $p\mid x_3x_4.$ By the construction, for any $p\mid x_3$ or $p\mid x_4$ we have that $p\nmid a,$ which is a contradiction.
\spm

\noindent{\bf Case 3.} Assume that $p\in\mathcal{P}_3.$ As in Case 2, $p\mid ax_2^2,$ so $p\mid x_2.$ But, for any $p\mid x_2$ we have that $p\nmid c,$ which is a contradiction.
\spm

\noindent{\bf Case 4.} Assume that $p\in\mathcal{P}_4.$ Since $p\mid a'$ and $p\mid x_4,$ we have that $p\mid ax_2^2,$ so $p\mid x_2,$ which is not possible since for any $p\mid x_2$ we have that $p\mid a.$ Hence, we obtain that $\gcd(a',d)=1.$

Observe that since $d$ is even, $b$ is also even. If $c$ is even, then $2\mid x_3.$ Thus, it follows that $a'=q(x_2,x_3x_4)=ax_2^2+bx_2x_3x_4+cx_3^2x_4^2\equiv ax_2^2\Mod{4}\equiv a\Mod{4}$ since $x_2$ is an odd number. If $c$ is odd, then $2\mid x_4$, and  so it follows that $a'=q(x_2,x_3x_4)=ax_2^2+bx_2x_3x_4+cx_3^2x_4^2\equiv ax_2^2\Mod{4}\equiv a\Mod{4}$ by the same reason. Hence, the first assertion follows.

For the second assertion, let $b=2b'$, for some integer $b'$, and observe that since $16 \mid (2b')^2-4ac$, we have $4\mid b'^2-ac$. Thus, we observe that if $c$ is even, then $4\mid 2b'$. If $c$ is odd, then $b'$ is odd. Hence, we similarly show that if $c$ is even, then $a'=q(x_2,x_3x_4)=ax_2^2+bx_2x_3x_4+cx_3^2x_4^2\equiv ax_2^2\Mod{8}\equiv a\Mod{8}$ since $4\mid b$, $2\mid x_3$ and $2\nmid x_2$. Also, if $c$ is odd, then 
$$
a'=q(x_2,x_3x_4)=ax_2^2+2b'x_2x_3 2 \frac{x_4}{2}+cx_3^2 4 \frac{x_4^2}{4}\equiv a + 4(b'x_2x_3\frac{x_4}{2}+cx_3^2  \frac{x_4^2}{4})\equiv a\Mod{8}
$$
since $b'x_2x_3\frac{x_4}{2}$ and $cx_3^2  \frac{x_4^2}{4}$ are odd. Thus, the second assertion follows.
\end{proof}

Given $\Omega$, $\Delta$ and $a$ under the hypotheses in (\ref{eq: construction_a}), we constructed a binary form $q_{\Delta,\Omega}$; cf.\ (\ref{eq: construction_q_Omega_Delta}), and an abelian surface $A = E_1\times E_2$ such that $q_{\Delta,\Omega} \sim q_{E_1,E_2}$ (cf.\ Proposition \ref{p: construction_A_th_1}, and the discussions above). Recall that $R_{3,4}(q_{\Delta,\Omega})\neq \varnothing$ by the hypothesis. 

In addition, suppose now that $\disc(q_{\Delta,\Omega}) = -2\Delta\Omega^2 \equiv 16 \Mod{32}$, and thus $q_{\Delta,\Omega}$ primitively represents two numbers $r_1'$ and $r_2'$ such that $r_1'\equiv 3 \Mod 8$ and $r_2'\equiv 7 \Mod 8$ by Lemma \ref{representstwonumbers}. Since $\Omega$ is odd, we have that $16\mid -2\Delta$. Thus, a primitive form $\frac{1}{\Omega}q_{\Delta,\Omega} = q_{\Delta}$ of discriminant $-2\Delta$ primitively represents $r_i/\Omega$, for $i=1,2$ such that $r_i/\Omega \equiv r_i'/\Omega \Mod{8}$ with $\gcd(r_i/\Omega, -2\Delta) =1$ by Lemma \ref{representsrelatively}. Since $\Omega$ is odd, we have that $r_i \equiv r_i'\Mod{8}$. Let us put $q_{E_1,E_2}(h_i):=r_i,$ for primitive isogeny $h_i\in\Hom(E_1,E_2),$ for $i=1,2.$ (Note that $r_i \in R(q_{\Delta,\Omega})=R(q_{E_1,E_2})$). Now, 
\begin{equation}
\label{basicpp}
    \theta_{r_i} \;=\; \textbf{D}(2,(r_i+1)/2,h_i) \;\in\; \mathcal{P}(A)^{\ev}, \textrm{ for } i=1,2
\end{equation} 
by Proposition \ref{eventhetaprop}.  Thus, we have constructed principally polarized abelian surfaces $(A,\theta_{r_i})$ such that 
\begin{equation}
\label{eq: notation_f'_Delta_Omega}
   f_{q_{\Delta,\Omega}}^{\theta_{r_i}}\;:=\; q_{(A,\theta_{r_i})} \text{ is an imprimitive ternary form for } i=1, 2,
\end{equation}
(by Corollary \ref{imprimitiveiffeven}). Let us state what we proved by this construction:

\begin{proposition}
\label{p: construction_A_th_2}
    For given positive odd $\Omega$ and even $\Delta$ integers, suppose that there is an integer $a$ satisfying the conditions in \emph{(\ref{eq: construction_a})}. In addition, suppose that $-2\Delta \Omega^2 \equiv 16 \Mod{32}$. Therefore, there is a product abelian surface $A = E_1 \times E_2$, where  $q_{E_1,E_2}$ is a binary quadratic form of discriminant $-2\Omega^2\Delta$ with content $\Omega$,  and there are two principal polarizations $\emph{\textbf{D}}(2, (r_i + 1)/2, h_i) \in \mathcal{P}(A)^{\ev},$ for some primitive isogenies $h_i \in \Hom(E_1,E_2)$,  for $i = 1, 2$. Also, we have that $\gcd(r_i/\Omega, 2\Delta) = 1$, $r_1 \equiv 3 \Mod{8}$ and $r_2 \equiv 7 \Mod{8}$.
\end{proposition}

In the following two important results, we use the technique of the use of assigned characters to show the genus equivalence. For this,
let $f_1$ and $f_2$ be primitive positive definite  ternary quadratic forms with the same basic (genus) invariants. As was mentioned in \S\ref{ternaryforms}, we have that $f_1$ and $f_2$ are genus-equivalent if and only if the following equations hold: 
\begin{equation}
\label{suffices_EXPLICIT}
     \chi(f_1) \;=\; \chi(f_2), \ \forall\chi\in X(f_1),\ \textrm{ and }\  \chi(F^B_{f_1})\;=\;\chi(F^B_{f_2}),\ \forall\chi\in X(F_{1}).
\end{equation}

By Proposition \ref{prop: existence_number_a} and Corollary \ref{cor: characters_for_ref_hum},  we prove the following  two useful results. They will be key for the main results in this article.

\begin{proposition}
\label{p: f+1.1+1.2_genus_eqv}
    Let $f$ be a positive definite ternary form satisfying conditions \emph{(\ref{classificationconditions1})} and \emph{(\ref{classificationconditions2})}. Let $\Omega=\Omega_{f/2}$ and $\Delta=\Delta_{f/2}$ be the basic invariants of $f/2$. Then $f$ is genus equivalent to $f_{q_{\Delta,\Omega}}^{\theta}$, where $\theta = \theta_a$, for some integer $a$; cf.\ \emph{(\ref{eq: construction_theta_for_Omega_Delta})}.
\end{proposition}

\begin{proof}
    By Proposition \ref{prop: existence_number_a}, there exists an odd integer $a$ satisfying the conditions in (\ref{eq: construction_a}) for $\Delta$, $\Omega$. By the constructions of the principally polarized abelian surfaces above (see Proposition \ref{p: construction_A_th_1}, and the discussions above), there is a product abelian surface $A=E_1 \times E_2$ such that $q_{E_1,E_2} \sim q_{\Delta, \Omega}$, where $q_{\Delta, \Omega}(1,0)= a\Omega$; cf.\ (\ref{eq: construction_q_Omega_Delta}), and there is $\theta :=\theta_a \in \mathcal{P}(A)^{\ev}$. Let us put $4f':=f^{\theta}_{q_{\Delta, \Omega}} = q_{(A,\theta)}$. By Corollary \ref{imprimitiveiffeven}, $f'$ is primitive.

    By the construction and (\ref{basicinvariantsDicksonsense}), we see that  $\Omega_{2f'} = \cont(q_{\Delta, \Omega})= \Omega$ and $\Delta_{2f'} = -\disc(q_{\Delta})/2 = \Delta$.  Therefore, it suffices to show that (\ref{suffices_EXPLICIT}) holds for $f/4$ and $f'$.

   We have by Proposition \ref{prop: existence_number_a} that $\chi(f/4) = 1$, for all $\chi\in X(f/4)$, and also, we have by Corollary \ref{cor: characters_for_ref_hum} that $\chi(f') = 1$ for all $\chi\in X(f')$. Thus, it follows that $\chi(f/4) = \chi(f')$, for all  $\chi\in X(f/4)$, and so the first equation of  (\ref{suffices_EXPLICIT}) holds for $f/4$ and $f'$.

Recall from Proposition \ref{F=F^B} that $F_{f/4}^B =F_{f/2}$ and $F_{f'}^B = F_{2f'}$.     We know that $\chi_p(a) = \chi_p(F_{f/2})$ for all odd primes $p\mid \Delta$ by Proposition \ref{prop: existence_number_a}. By the definition of $\theta$ (see (\ref{eq: construction_theta_for_Omega_Delta})), it follows from Proposition \ref{invariants}  that $F_{2f'}$ represents $a$. Since $\gcd(a,\Delta)=1$ by (\ref{eq: construction_a}), it follows that
$$
\chi_p(F_{2f'}) \;=\; \chi_p(a)  \;=\; \chi_p(F_{f/2}), \ \text{ for all odd primes } \ p\mid \Delta.
$$

Since $\Delta$ is even, we have that $\Delta \equiv 0, 2 \Mod{4}$. If $\Delta\equiv 2 \Mod{4}$, then $X_s(F_{2f'})=\{\chi_{-4}\}$ by (\ref{eq: assigned_character_F_f_Smith_sense}), and it follows by (\ref{Smithcharacter}) that $\chi_{-4}(F_{2f'})=\chi_{-4}(\Omega)= \chi_{-4}(F_{f/2})$. 
If $\Delta\equiv 0 \Mod{4}$, then $X_s(F_{2f'})=\{\chi_{-4}, \chi_8\}$ by (\ref{eq: assigned_character_F_f_Smith_sense}), and we have again by (\ref{Smithcharacter}) that $\chi_{-4}(F_{2f'})=\chi_{-4}(F_{f/2})$. Moreover, we see from Proposition \ref{prop: existence_number_a} that $\chi_8(a)=\chi_8(F_{f/2})$. Since $F_{2f'}$ represents an odd integer $a$, we also have that $\chi_8(a)=\chi_8(F_{2f'})$, and hence $\chi_8(F_{2f'}) = \chi_8(a)  = \chi_8(F_{f/2})$. Therefore,  we have shown that  (\ref{suffices_EXPLICIT}) holds for $f/4$ and $f'$ in both cases, and this implies that $f/4$ and $f'$ are genus equivalent, which proves the assertion.
\end{proof}


\begin{theorem}
\label{cor: Theta_lies_genus}
    Let $A=E_1\times E_2$ be a product surface, where $E_1 \sim E_2$ are CM elliptic curves, and let $\mathcal{P}(A)^{\ev}\neq\varnothing.$ Let $d:=-\disc(q_{E_1,E_2})$. Then we have that
    \begin{enumerate}
\item[\emph{(i)}] if $d \not\equiv 16 \Mod{32},$ then $\Theta(A)^{\ev}$ lies in a single genus,
\item[\emph{(ii)}] if $d \equiv 16 \Mod{32},$ then $\Theta(A)^{\ev}$ lies in exactly two genera.
\end{enumerate}
\end{theorem}

\begin{proof}
    Let us take two arbitrary principal polarizations $\theta,\theta'\in\mathcal{P}(A)^{\ev},$ and put $f_1:=\frac{1}{2}q_{(A,\theta)}$ and $f_2:=\frac{1}{2}q_{(A,\theta')}$. We have by Proposition \ref{invariants} that the $f_i$ are improperly primitive, so let $F_i:=F_{f_i}$ be the reciprocal of $f_i$ for $i=1,2.$   Also, we know by Proposition \ref{F=F^B} that the reciprocals $F_i=F_{f_i/2}^B$ are properly primitive. 
    
    By Corollary \ref{cor: characters_for_ref_hum} (and by (\ref{I1Omega})), $f_1$ and $f_2$ have the same basic invariants $\Omega$ and $\Delta$.  We also know by this corollary  that  
\begin{equation}
\label{eq: the_claim_in_proof}
    \chi(f_1)=\chi(f_2),\forall\chi\in X(f_1), \textrm{ and }  \chi(F_{1})=\chi(F_{2}),\forall\chi\in X(F_{1})\setminus\{\chi_8, \chi_{-4}\chi_8\}.
\end{equation}

Let $t:=\cont(q_{E_1,E_2})$ and $d':=d/t^2.$  We have by Proposition \ref{invariants} that $\Omega=t$ is odd and $\Delta=d'/2$ is even.

(i)  If $d=t^2d'\not\equiv16 \Mod{32}$, then $d'\not\equiv 16 \Mod{32}$ since $t$ is odd. Also, since $\Delta=d'/2$ is even, we have that $\Delta\equiv 0,2\Mod4$. If $\Delta\equiv 2\Mod4$, then we conclude from (\ref{eq: the_claim_in_proof})  that (\ref{suffices_EXPLICIT}) holds since $\chi_{-4}$ is the only supplementary assigned character of $F_i$ in this case by (\ref{eq: assigned_character_F_f_Smith_sense}). This means that $f_1$ and  $f_2$ are genus-equivalent, which proves the assertion (i) in this case. 

Assume now that $\Delta\equiv0\Mod4,$ so $2\Delta=d'\equiv 0 \Mod{8}$. Then since   $d'\not\equiv 16 \Mod{32}$, we have that  $d'\equiv 0, 8, 24 \Mod{32}$. To prove the assertion (i) in this situation, we consider two cases below, namely $d'\equiv 0, 24 \Mod{32}$ or $d'\equiv 8 \Mod{32}$. Before starting the analysis of these cases, note that $\chi_{-4}, \chi_8\in X_s(F_1)$ in this case by (\ref{eq: assigned_character_F_f_Smith_sense}).  Hence, it is clear that $X_s(q')\subseteq X_s(F_1)$, where  $q':=\frac{1}{t}q_{E_1,E_2}.$

 Let us put $\theta=:\textbf{D}(n_1,m_1,k_1h_1)$ and $\theta'=:\textbf{D}(n_2,m_2,k_2h_2)$,  where the $h_i$ are primitive elements in $\Hom(E_1,E_2)$, and let $\delta_i:=q'(h_i)$ for $i=1, 2$. Since $\theta, \theta' \in\mathcal{P}(A)^{\ev}$, we have that $n_im_i-q_{E_1,E_2}(k_ih_i)=1$ by Corollary 25 of \cite{MJ}, and $n_im_i \equiv 0 \Mod{4}$ by (\ref{eventheta}), and so $q_{E_1,E_2}(k_ih_i)=k_i^2q_{E_1,E_2}(h_i)\equiv 3 \Mod{4}$. Thus, since $q_{E_1,E_2}(h_i)$ is odd, and since $q_{E_1,E_2}(h_i)=t\delta_i$, it follows that the $\delta_i$ are odd. 
\spm

\noindent{\bf Case 1.} Assume that $d'\equiv 0 \text{ or } 24 \Mod{32}.$
\sps

\noindent By \cite[p.~55]{davidcox}, $\chi_8\in X(q')$ in these cases, where $X(q')$ denotes the set of the assigned characters of $q'.$ Since $\delta_i$ is odd and $q'\rightarrow\delta_i,$ we have that $\chi_8(q')=\chi_8(\delta_i).$ Moreover, since $\frac{1}{t} q_{E_1,E_2}(h_i)=q'(h_i)=\delta_i,$  we obtain by Proposition \ref{invariants}  that $\delta_i$ is represented by the reciprocal $F_i$ for $i=1,2,$ so $\chi_8(F_i)=\chi_8(\delta_i)=\chi_8(q').$ In particular, $\chi_8(F_1)=\chi_8(F_2).$ Hence, it follows by (\ref{eq: the_claim_in_proof}) that (\ref{suffices_EXPLICIT}) holds, so $f_1=\frac{1}{2}q_{(A,\theta)}$ and $f_2=\frac{1}{2}q_{(A,\theta')}$ are genus-equivalent, and so $q_{(A,\theta)}$ and $q_{(A,\theta')}$ are genus-equivalent. Since $\theta_1$ and $\theta_2$ are arbitrary elements in $\mathcal{P}(A)^{\ev},$ the assertion (i) follows in this case.
\spm

\noindent{\bf Case 2.}  Assume that $d'\equiv 8 \Mod{32}.$
\sps

\noindent By \cite[p.~55]{davidcox}, $\chi_{-4}\chi_8\in X(q')$ in this case. As in Case 1, we obtain that $\chi_{-4}\chi_8(q')=\chi_{-4}\chi_8(\delta_i),$  for $i=1,2,$ and  $\chi_{-4}\chi_8(F_1)=\chi_{-4}\chi_8(F_2).$ Since  $\chi_{-4}(F_1)=\chi_{-4}(F_2)$ by (\ref{eq: the_claim_in_proof}),
it follows that $\chi_8(F_1)=\chi_8(F_2)$ again, so the assertion (i) follows in a similar way as in Case 1. Thus, the assertion (i) holds in all cases.

(ii) Assume now that $d \equiv 16 \Mod{32}$. By the construction of the principal polarizations in (\ref{basicpp}), we can see that there are two principal polarizations $\theta_{r_1}$ and $\theta
_{r_2} \in \mathcal{P}(A)^{\ev}$ such that $\gcd(r_i,2\Delta)=1$ and $r_1\equiv r_2+4 \Mod{8}$ (see Proposition \ref{p: construction_A_th_2}).

  We first observe that the improperly primitive forms $f_1':=\frac{1}{2}q_{(A,\theta_{r_1})}$ and $f_2':=\frac{1}{2} q_{(A,\theta_{r_2})}$ are not genus-equivalent. Indeed,  since $t$ is odd, we have that  $d' \equiv 16 \Mod{32}$. Then, since $\Delta_{f_i'}=d'/2$ by (\ref{basicinvariantsDicksonsense}), and since $d'\equiv 16\Mod{32},$ we have that $4\mid \Delta_{f_i'}$, and hence $\chi_8\in X(F_i'),$ where $F_i'$ is the reciprocal form of the $f_i'$'s, by (\ref{eq: assigned_character_F_f_Smith_sense}). By Proposition \ref{invariants}, we have that $F_i'\rightarrow r_i/t.$ Since $t$ is odd by what was mentioned above (or by \ref{todd4dividesd'}) and since $t\equiv t^{-1} \Mod{8}$, we see that $\chi_8(r_i/t)=\chi_8(r_i)\chi_8(t)$. Thus, since  $r_1\equiv r_2+4 \Mod 8$, we get that
\begin{align}
\label{eq: in_proof_not_genus_equv}
    \chi_8(F_1') \;=\; \chi_8(r_1)\chi_8(t) \;=\; -\chi_8(r_2)\chi_8(t) \;=\; -\chi_8(r_2/t) \;=\; - \chi_8(F_2').
\end{align}
Therefore, it follows that  $\chi_8(F_1')\neq\chi_8(F_2'),$ which means that (\ref{suffices_EXPLICIT}) does not hold, and so $q_{(A,\theta_{r_1})}\text{ and } q_{(A,\theta_{r_2})}$ are not genus-equivalent.


Moreover, for any $\theta''\in\mathcal{P}(A)^{\ev}$, $\chi_8(F_{\frac{1}{2}q_{(A,\theta'')}}) = \chi_8(F_1')$ or $\chi_8(F_2')$. Therefore, we see from (\ref{eq: the_claim_in_proof}) and (\ref{suffices_EXPLICIT}) that either $\frac{1}{2}q_{(A,\theta'')}$ and $f_1'$ or $\frac{1}{2}q_{(A,\theta'')}$ and $ f_2'$ are genus equivalent. Hence, it follows that any $q_{(A,\theta'')}\in\Theta(A)^{\ev}$ lies in either $\gen(2f_1')$ or $\gen(2f_2')$, so the  assertion (ii) follows.
\end{proof}

\section{Classification of the refined Humbert invariant}
\label{section: proofs_main results}

In this section, the aim is to prove the two main results, namely Theorems \ref{theorem1} and \ref{maintheorem} in this article. 
Given a positive definite ternary form $f$ satisfying conditions (\ref{classificationconditions1}) and  (\ref{classificationconditions2}), we  show in this section that there exists a principally polarized abelian surface $(A,\theta)$ such that $f\sim q_{(A,\theta)}.$ This proves (the hard part of) Theorem \ref{maintheorem}. To this end, we want to apply Theorem 34 of \cite{dicksonsbook}. For this, we need to find a suitable binary form which is primitively represented by a ternary refined Humbert invariant and $f$. We see that this binary form can be of a form $q_s,$ for some $s\in P(n)^{\ev}$, for some integer $n$; cf.\ Corollary \ref{q_sphi_p} below. This observation requires considerable work regarding the theory of binary forms represented by ternary forms. 

The following lemma will be useful in the proof of Proposition \ref{thereisformrepresentsp}. 

\begin{lemma} 
\label{thesetnonempty}
Let $A=E_1\times E_2$ be a product surface, where $E_1\sim E_2$ are CM elliptic curves such that $\mathcal{P}(A)^{\ev}\neq\varnothing.$  Let $tq' = q_{E_1,E_2} $, for some primitive form $q'$. Then there exists an integer $r$ relatively prime to $\disc(q')$ such that both $q'$ and $F_f$  represent $r$, where $2f = q_{(A,\theta)}$, for some $\theta \in \mathcal{P}(A)^{\ev}$.
\end{lemma}

\begin{proof}
 We first claim that there is a number $r$ represented by $q'$ such that $tr\equiv 3 \Mod{4}$ and $\gcd(r,\disc(q')) = 1$. Since $\mathcal{P}(A)^{\ev}\neq\varnothing$, we have that $t$ is odd, $4\mid d':=\disc(q')$ and $R_{3,4}(q_{E_1,E_2})\neq \varnothing$ by Proposition \ref{eventhetaclassification}. Since $R_{3,4}(q_{E_1,E_2})\neq \varnothing$, we may clearly suppose that $q_{E_1,E_2} \rightarrow r'$ with $r'\equiv 3 \Mod{4}$, and so $q' \rightarrow r'':=r'/t$. Hence, $q' \rightarrow r$ with $r\equiv r'' \Mod{4}$ and $\gcd(r,d') = 1$ by Lemma \ref{representsrelatively}. Thus, $tr \equiv tr''\equiv r' \equiv 3\Mod{4}$, and this proves the claim. 
 
 Therefore, there is  a primitive isogeny $h\in\Hom(E_1,E_2)$ such that $q_{E_1,E_2}(h)=tr\equiv 3 \Mod{4}.$   Now, consider $\theta=\textbf{D}(2,(tr+1)/2,h)\in\mathcal{P}(A)^{\ev}$ (cf.\ Proposition \ref{eventhetaprop}). Let $f:=\frac{1}{2}q_{(A,\theta)},$ and let $F_{f}$ be its reciprocal. Then we know that $F_f$ represents $r$ by Proposition \ref{invariants}. 
 Hence $r\in R(q')\cap R(F_{f})$ with $\gcd(r,d')=1,$ which proves the assertion. 
 \end{proof}

 The following result finds a binary form  $q\in\gen(\frac{1}{t}q_{(E_1, E_2)})$ such that $q$ and $F_f$ represent the same prime number when we are in the situation of Proposition \ref{invariants}.

\begin{proposition}
\label{thereisformrepresentsp} 
Let $A=E_1\times E_2$, $d \text{ and } t$ be as in \emph{Proposition \ref{invariants}}. Let $\theta \in \mathcal{P}(A)^{\ev}$, and let $f\in\gen(\frac{1}{2}q_{(A,\theta)}).$ Then there exists a prime number $p$ represented by the reciprocal $F_f$ of $f$ with $p\nmid d.$ Moreover, for such a number $p$  there exists a binary form  $\tilde{q}'\in\gen(q')$ which also represents $p,$ where $q':=\frac{1}{t}q_{E_1,E_2}.$
\end{proposition}

\begin{proof}
    First, we know that the content and the basic invariants $I_1$ and $I_2$ are genus invariants. Thus, since we have from Proposition \ref{invariants} that $f$ is improperly primitive, it follows that $F_f$ is properly primitive by Proposition \ref{F=F^B}, and so it represents infinitely many prime numbers by Corollary 12 of \cite{refhum}. In particular, there is a prime number $p\nmid d$ represented by $F_f$. This proves the first assertion.

    Let $d' = d/t^2$ (as in Proposition \ref{invariants}), and let $X(q')$ be the set of assigned characters of $q'$ as defined in \cite[p.~55]{davidcox}.  Since $q'$ is a primitive form,  the second assertion follows from Theorem 2.26 and Lemma 3.20 of \cite{davidcox} (or from Theorem 4.4 of \cite{Duncanbinaryforms}), once we have shown that 
\begin{align}
\label{chi(p)=chi(q)}
    \chi(p)\;=\;\chi(q'), \text{ for all } \chi\in X(q').
\end{align} 

Let $X_o(q'):=\{\chi_{\ell}: \ell\mid d', \ell\neq2 \text{ prime}\},$ and let $X_s(q')$ denote the set of the supplementary assigned characters of $q'$, so $X(q')=X_o(q')\cup X_s(q').$ By Proposition \ref{F=F^B}, we have that $F^B_{f/2}=F_f.$ Let $X_o(F_f):=\{\chi_{\ell}: \ell\mid \Delta_f, \ell\neq2 \text{ prime}\},$ and let $X_s(F_f)$ denote the supplementary assigned characters of $F_f.$ Since $\Delta_f=d'/2$ by Eq.~(\ref{basicinvariantsDicksonsense}), we obtain that $X_o(q')=X_o(F_f).$ 

We have from Lemma \ref{thesetnonempty} (and its proof) that there exists $\theta^{\prime} = \textbf{D}(n, m, k h)\in P(A)^{\ev}$, for a primitive $h \in \Hom(E_1, E_2)$ such that $q^{\prime}(h)=a^{\prime}$ is relatively prime to $d^{\prime}$, and $a^{\prime}$ is represented by $F_{f^{\prime}}$, where $f^{\prime}=\frac{1}{2}q_{(A, \theta^{\prime})}$. 
Since $\Delta_f=\Delta_{f'}$ by (\ref{basicinvariantsDicksonsense}), we have that $X_o(F_{f'})=X_o(F_{f})$. Therefore, we obtain that $\chi(q')=\chi(a')=\chi(F_{f'})$, for all $\chi\in X_o(q')=X_o(F_{f'})$. Moreover, since $\chi(F_{f})=\chi(F_{f'})$ by Corollary \ref{cor: characters_for_ref_hum}, we obtain that 
\begin{equation}
\label{eq2: characters_of_reciprocals_in_proof}
    \chi(p) \;=\; \chi(F_f) \;=\; \chi(F_{f'}) \;=\; \chi(a') \;=\; \chi(q'), \quad \forall \chi \in X_o(F_{f}).
\end{equation}

In order to prove that $\chi(q')=\chi(p)$, for all  $\chi\in X_s(q'),$ we  use the supplementary assigned characters $X_s(F_f)$  of $F_f$. Note that since $\Delta_f$ is even, we know that $\chi_{-4}\in X_s(F_f)$ by (\ref{eq: assigned_character_F_f_Smith_sense}).

Observe now that if $\chi_{-4}\in X_s(q')\cap X_s(F_f)$, then $\chi_{-4}(p)=\chi_{-4}(q')$ because (\ref{cor: eq_characters}) implies that
\begin{equation*}
    \chi_{-4}(p) \;=\; \chi_{-4}(F_f) \;=\; \chi_{-4}(F_{f'}) \;=\; \chi_{-4}(a') \;=\; \chi_{-4}(q').
\end{equation*}

If $d^{\prime} \equiv 16 \Mod{32}$, then $X_s(q^{\prime})=\{\chi_{-4}\}$ (by \cite[p.~55]{davidcox}) and $\chi_{-4} \in X_s(F_f)$, and  so 
our observation, together with (\ref{eq2: characters_of_reciprocals_in_proof}), proves (\ref{chi(p)=chi(q)}), i.e, the second assertion follows in this case.

We henceforth suppose that $d^{\prime} \not\equiv 16 \Mod{32}$. If $d^{\prime} \equiv 4 \Mod{8}$, then we see that $X_s(q')=\varnothing$ or $\{\chi_{-4}\}$ when $d' / 4 \equiv 3 \Mod{4}$ or $1\Mod{4}$ respectively by \cite[p.~55]{davidcox}.  When $X_s(q')=\varnothing$, (\ref{chi(p)=chi(q)}) already holds. When $X_s(q')=\{\chi_{-4}\}$,  since $\chi_{-4} \in X_s(F_f)$,  (\ref{chi(p)=chi(q)}) follows again by the observation (together with (\ref{eq2: characters_of_reciprocals_in_proof})).

Suppose now that $d^{\prime} \equiv 0\Mod{8}$ (and still $d'\not\equiv 16 \Mod{32}$).  We have from (\ref{eq: assigned_character_F_f_Smith_sense})
that $X_s(q') \subset X_s(F_f)$. Since $f$ and $f'$ are genus equivalent by Theorem \ref{cor: Theta_lies_genus}, we have that  $\chi(F_f) = \chi(F_{f'})$ for all $\chi\in X_s(F_f)$. Thus, it follows as in (\ref{eq2: characters_of_reciprocals_in_proof}) that 
$\chi(p) = \chi(F_f) = \chi(F_{f'}) = \chi(a') = \chi(q')$, for all $\chi \in X_s(q')\subseteq X_s(F_f)$. Hence, (\ref{chi(p)=chi(q)}) holds again by (\ref{eq2: characters_of_reciprocals_in_proof}), and so the assertion follows.
\end{proof}

\begin{proposition} 
\label{frepresentsphi}
Let $A=E_1\times E_2$, $d$ and $t$ be as in \emph{Proposition \ref{invariants}}, and let $\theta\in\mathcal{P}(A)^{\ev}.$ For any $f\in\gen(\frac{1}{2}q_{(A,\theta)}),$ there is a prime number $p\nmid d$ represented by $F_f,$ and there exists a binary form $\phi_p$ of discriminant $-4\Omega_{f}p=-4tp$ such that $f\rightarrow \phi_p.$ Moreover, $\cont(\phi_p)=2.$
\end{proposition}

\begin{proof} The first assertion follows from Proposition \ref{thereisformrepresentsp}, and so there is a prime number $p$ with $p\nmid d,$ which is represented by $F_f.$   
We have by Eq.~(\ref{basicinvariantsDicksonsense}) that $\Omega_f=t$ and $\Delta_{f}=d'/2$.
 For any such $p$ we have by Theorem 38 of \cite{dicksonsbook} that there is a binary form $\phi_p$ with $\det^{\text{D}}(\phi_p)=-\Omega_f p=-tp$ such that $f\rightarrow \phi_p.$ Hence $\disc(\phi_p)=-4tp$ (cf.\ Eq.~(\ref{dicksondeterminant})) with $p\nmid d,$ which proves the second assertion. Moreover, we see from Theorem 37 of \cite{dicksonsbook} that $\phi_p$ is a properly or improperly primitive form since $\gcd(p,\Omega_f\Delta_f)=\gcd(p,d/2)=1.$ Also, note that $\cont(f)=2$ by Proposition \ref{invariants}, so $2\mid\cont(\phi_p).$ Therefore, $\phi_p$ is improperly primitive, and hence, $\phi_p/2$ is primitive (in Watson's sense), i.e., $\cont(\phi_p)=2,$ which proves the last assertion.  \end{proof}

The following result leads to a useful fact related to the forms $q_s;$ cf.\ Corollary \ref{q_sphi_p} below.

\begin{theorem}
\label{phiprincipalgenus} Let $A,t,d'$ and $\theta$ be as in \emph{Proposition \ref{invariants}}.
Let $f\in\gen(\frac{1}{2}q_{(A,\theta)}).$ If $f$ primitively represents a binary form $\phi_p$ with $\cont(\phi_p)=2$ and $\disc(\phi_p)=-4tp,$ for some prime $p\nmid td',$ then $\phi':=\frac{1}{2}\phi_p\in\gen(1_{-tp}),$ i.e., $\phi'$ lies in the principal genus of discriminant $-tp.$ 
\end{theorem}

\begin{proof}
Since $\cont(\phi_p)=2$ and $\disc(\phi_p)=-4tp$, we see that $\phi'$ is a primitive form of discriminant $-tp$. 
By Lemma 3.20 of \cite{davidcox}, it suffices to show that $\chi(\phi')=1$ for all assigned characters $\chi\in X(\phi').$  Since $tp$ is odd (cf.\ Eq.~(\ref{todd4dividesd'})), we obtain that the assigned characters of $\phi'$ are  the $\chi_{\ell}$'s for odd primes $\ell\mid tp;$ cf.\ \cite[p.~55]{davidcox}. Moreover, the assigned characters of $f/2$ are the $\chi_{\ell}$'s with primes $\ell\mid t$ since $I_1(f/2)=-t$ by Eq.~(\ref{I1I2basicinvariants}). Thus, we have from  (\ref{cor: eq_characters}) that
\begin{align}
\label{allchsare1}
    \chi_{\ell}(f/2)\;=\;1, \textrm{ for any prime  } \ell\mid t.
\end{align}

Let $q:=q_{E_1,E_2}$. There exists an integer $r$ represented by $\phi'$ with $\gcd(r,td')=1$ by Lemma 2.25 of \cite{davidcox}. Since $f/2\rightarrow\phi',$  $r$ is also represented by $f/2.$ By equation (\ref{allchsare1}), we therefore obtain that \begin{equation}
 \label{Claimpart}
    1\;=\; \chi_{\ell}(f/2) \;=\; \chi_{\ell}(r) \;=\; \chi_{\ell}(\phi'), \ \text{ for all primes  } \ \ell\mid t.
\end{equation} 
So it remains to show that $\chi_p(\phi')=1.$ Since $\disc(\phi')=-tp$ is an odd number, there are infinitely many primes $p'\equiv 1\Mod{4}$ represented by $\phi'$ by Corollary 4 of \cite{refhum}, so there is a prime $p'$ with $p'\nmid tp$ and  $p'\equiv 1\Mod{4}$ such that $\phi'\rightarrow p'.$  Then by Lemma 2.5 of \cite{davidcox}, we have that $\left(\frac{-tp}{p'}\right)=1.$ Since $p' \in R(\phi') \subset R(f/2)$, so we get that 
\begin{equation} 
\label{allprimes}
    \chi_{\ell}(p') \;=\; \chi_{\ell}(f/2)\;=\;1 \text{ for all primes } \ell\mid t
\end{equation} by Eq.~(\ref{Claimpart}). By the Quadratic Reciprocity Law, we have that 
\begin{equation}
\label{reciprocity}
\left(\frac{p'}{\ell}\right) \;=\; \left(\frac{\ell}{p'}\right), \textrm{ for all } \ell, \text{ and also }  \left(\frac{p'}{p}\right)=\left(\frac{p}{p'}\right) 
\end{equation} 
because $p'\equiv 1\Mod{4}.$ 

Let us define the set $S_t:=\{\ell \text{ prime }: v_{\ell}(t) \text{ is odd }\},$ i.e., $S_t$ is the set of prime numbers which appear with an odd power in the factorization of $t.$ We therefore get that 
\begin{align*}
1&= \left(\frac{-tp}{p'}\right) =  \left(\frac{p}{p'}\right)\prod_{\ell\in S_t}\left(\frac{\ell}{p'}\right)\stackrel{\text{(\ref{reciprocity})}}{=}\left(\frac{p'}{p}\right)\prod_{\ell\in S_t}\left(\frac{p'}{\ell}\right)\stackrel{\text{(\ref{allprimes})}}{=}\left(\frac{p'}{p}\right)=\chi_p(p') \;=\;\chi_p(\phi').&
\end{align*}
Hence we have proved that $\chi_p(\phi')=1.$ This, together with Eq.~(\ref{Claimpart}), proves that $\phi'$ lies in the principal genus of discriminant $-tp.$ \end{proof}

\begin{corollary}
\label{q_sphi_p}
Assume that we are in the situation of \emph{Theorem \ref{phiprincipalgenus}.} Then there exists $s\in P(tp)^{\ev}$ such that $q_s\sim2\phi_p.$
\end{corollary}

\begin{proof}
By the hypothesis and by the result of Theorem \ref{phiprincipalgenus}, $\phi_p$  is  a binary quadratic form  with $\cont(2\phi_p)=4$ and $\disc(2\phi_p)=-16tp$ such that  $\frac{1}{2}\phi_p$ lies in the principal genus. Since the discriminant of $\frac{1}{2}\phi_p$ is $-tp$ is odd (cf.\ the proof of Theorem \ref{phiprincipalgenus}), we have that $-tp\equiv1\Mod{4}$.  We therefore get that $2\phi_p\sim q_s,$ for some $s\in P(tp)^{\ev}$ by Proposition \ref{classificationoftypedelta}. Hence, the assertion follows.
\end{proof}

Recall from \cite{milne1986jacobian} that if $C/K$ is a curve of genus $2$, then its Jacobian $J_C$ is an abelian surface and there is a divisor $\theta_C$ on $J_C$, called the \textit{theta-divisor} such that $\theta_C$ is a principal polarization in $\NS(J_C)$, and $\theta_C\simeq C$. We then write $q_C := q_{(J_C,\theta_C)}$ for its associated refined Humbert invariant as was mentioned in the introduction.  One key property of this invariant $q_{(A,\theta)}$ is determining whether $(A,\theta)$ is a Jacobian. By Proposition 6 of \cite{MJ}, we have that 
\begin{equation} 
\label{thetaisirreducible}
(A,\theta) \simeq (J_C,\theta_C), \mbox{ for some curve }C/K\; \Leftrightarrow\; 
q_{(A,\theta)}(D)\neq 1, \forall D\in \NS(A,\theta).
 \end{equation}

We are now ready to prove the first main result.
\spm

\noindent\textbf{Proof of Theorem \ref{theorem1}.} Since $f_1$ is an imprimitive ternary form, which is equivalent to a form $q_C$, there are isogenous CM elliptic curves $E_1,E_2$ such that $J_C=E_1\times E_2$ is a product surface by Proposition \ref{imprimitiveform} and $\theta:=\theta_C\in\mathcal{P}(J_C)^{\ev}$  by Corollary \ref{imprimitiveiffeven}. Let us put (as we did in the previous results) 
$$
 q \;:=\; q_{E_1,E_2}, \ d \;:=\; -\disc(q), \ t \;:=\; \cont(q), \textrm{ and } d' \;:=\; d/t^2.
$$ 
Suppose that $f_2\in\gen(f_1).$ We first prove the hard part that 
\begin{align}
\label{hardpart}
    f_2 \;\sim\; q_{\tilde{C}}, \  \textrm{ for some curve } \ \tilde{C} \text{ of genus 2},
\end{align} where $J_{\tilde{C}}$ is a product surface and $\theta_{\tilde{C}}\in\mathcal{P}(J_{\tilde{C}})^{\ev}.$ We apply Proposition \ref{frepresentsphi} to $f_2/2\in\gen(f_1/2)=\gen(\frac{1}{2}q_{C}).$ Thus, there is a prime number $p\nmid td'$ represented by $F_{f_2/2},$ and there is a binary form $\phi_p$ of discriminant $-4\Omega_{f_2/2}p=-4tp$ such that $f_2/2\rightarrow \phi_p,$ and $\cont(\phi_p)=2.$ By Corollary \ref{q_sphi_p},  
 we see that $2\phi_p$ has \textit{type} $tp,$ i.e., $2\phi_p\sim q_s,$ for some $s=(n,m,k)\in P(tp)^{\ev}.$ Since $f_2/2\rightarrow \phi_p$ and $\phi_p\sim q_s/2,$ we have that $f_2/2\rightarrow q_s/2$ by Theorem 28 of \cite{dicksonsbook}.

By Proposition \ref{thereisformrepresentsp}, there exists a binary form $\tilde{q}'\in\gen(q')$ with $\tilde{q}'\rightarrow p,$ where $q':=\frac{1}{t}q.$ We also have that there exist two elliptic curves $E_1',E_2'$ with $q_{E_1',E_2'}\sim \tilde{q}:=t\tilde{q}'$ by Lemma \ref{lemma26}. Since $\tilde{q}'$ primitively represents $p,$ $\tilde{q}$ primitively represents $tp,$ and thus, there exists a primitive $h\in\Hom(E_1',E_2')$ such that $q_{E_1',E_2'}(h)=tp.$

Let us put $A':=E_1'\times E_2',$ and let $\theta':=D_{s,h}:=\textbf{D}(n,m,kh)\in\mathcal{P}(A');$ cf.\ Corollary 25 of \cite{MJ}. Thus, we have that $q_{(A',\theta')}\rightarrow q_s$ by Proposition 29 of \cite{MJ}.

Since $s\in P(tp)^{\ev}$,  both $n$ and $m$ are even by definition. In addition, we have that $4 \mid d'$  by (\ref{todd4dividesd'}). Hence, $4\mid -t^2d'=\disc(q_{E_1',E_2'})$ since $q_{E_1',E_2'} \sim \tilde{q}$ and $\tilde{q} \in \gen(q)$. Therefore, $\theta'\in\mathcal{P}(A')^{\ev}$ by Eq.~(\ref{eventheta}). This implies that $q_{(A',\theta')}$ does not represent 1 because it is an imprimitive form by Corollary \ref{imprimitiveiffeven}, and so we have that $(A',\theta')\simeq (J_{\tilde{C}}, \theta_{\tilde{C}})$, for some curve $\tilde{C}$ of genus 2 by (\ref{thetaisirreducible}). Thus, if we let $f':=\frac{1}{2}q_{\tilde{C}}$, then it follows from Eq.~(\ref{basicinvariantsDicksonsense}) that $\Omega_{f'}=t$ and $\Delta_{f'}=d'/2$ since $\cont(q_{E_1',E_2'})=\cont(\tilde{q})=t$ and $\disc(q_{E_1',E_2'})=-t^2d'.$  Moreover, we know that  the basic invariants are genus invariants, and thus $f_1$ and $f_2$ have the same basic invariants.  Thus, since  $f_1\sim q_{C},$ we see that $\Omega_{f_2/2}=t$ and $\Delta_{f_2/2}=d'/2$ by Eq.~(\ref{basicinvariantsDicksonsense}) again, and so $f'$ and $f_2/2$ have the same basic invariants $\Omega:=\Omega_{f'}=\Omega_{f_2/2}$ and $\Delta:=\Delta_{f'}=\Delta_{f_2/2}$.

Since $\frac{1}{2}f_2$ and $f'$ are both improperly primitive (cf.\ Proposition \ref{invariants}) and have the same invariants $\Omega=t$ and $\Delta=d'/2,$ and since both primitively represent the binary form $\frac{1}{2}q_s$ whose discriminant is $-4tp,$ we can apply Theorem 34 of \cite{dicksonsbook}. Note that the determinant of $\frac{1}{2}q_s$ is $\det^{\text{D}}(\frac{1}{2}q_s)=-tp;$ cf.\ Eq.~(\ref{dicksondeterminant}). Therefore, we get from this theorem that $\frac{1}{2}f_2 \sim \frac{1}{2}q_{\tilde{C}} = f'$,  so $f_2\sim q_{\tilde{C}}$, which verifies Eq.~(\ref{hardpart}). 

To complete the proof, we have to show that $f_2 \sim q_{(J_C,\tilde{\theta})},$ for some $\tilde{\theta}\in\mathcal{P}(J_C)^{\ev}.$ To this end, note that $\tilde{q}\in\gen(q)$, and thus, it follows from Corollary 30 of \cite{refhum} that there exists a principal polarization $\tilde{\theta}\in\mathcal{P}(J_C)$ such that $ q_{(J_C,\tilde{\theta})}\sim q_{\tilde{C}}$.   Since $q_{\tilde{C}}$ is imprimitive, $q_{(J_C,\tilde{\theta})}$ is also imprimitive, and  so $\tilde{\theta}\in\mathcal{P}(J_C)^{\ev}$ by Corollary \ref{imprimitiveiffeven}. Hence, $f_2\sim q_{\tilde{C}} \sim q_{(J_C,\tilde{\theta})},$ for some $\tilde{\theta}\in\mathcal{P}(J_C)^{\ev},$ which proves the assertion because $\tilde{\theta}$ is isomorphic to a genus 2 curve $C'$ on $J_C$. 
\spm

Now, we are in a position to prove Theorem \ref{maintheorem}. 

\sps
\noindent\textbf{Proof of Theorem \ref{maintheorem}}.
($\Leftarrow$) By Proposition \ref{p: f+1.1+1.2_genus_eqv}, it follows that $f$ is genus equivalent to an imprimitive form $q_{(A,\theta)}$, for some $(A,\theta)$. Since $q_{(A,\theta)}$ is imprimitive, it follows from (\ref{thetaisirreducible}) that $q_{(A,\theta)} \sim q_C$, for some curve $C$ of genus 2. Thus, the assertion follows from  Theorem \ref{theorem1}. 

($\Rightarrow$) Assume that an imprimitive ternary form $f$ is equivalent to some refined Humbert invariant $q_{C}$. We first see that condition (\ref{classificationconditions2}) holds by Proposition \ref{necessaryconditionforii}. Secondly, it follows from Proposition \ref{imprimitiveform} and Corollary \ref{imprimitiveiffeven} that $J_C = E_1 \times E_2$ is a product surface, where $E_1\sim E_2$ are CM elliptic curves, and $\theta_C\in\mathcal{P}(J_C)^{\ev}.$ Thus, condition (\ref{classificationconditions1}) also holds by Proposition \ref{invariants}, and so the assertion follows.

Moreover, if this is the case, then as in the proof of Proposition \ref{imprimitiveform}, we have that $J_C\sim E\times E$, for some CM elliptic curve $E$.



\section{Applications on intersections of Humbert surfaces}
\label{section: applications}

In order to determine how the curves of genus $2$ whose Jacobians are isomorphic to some product abelian surface are distributed in the moduli space $\mathcal{M}_2(K)$  of genus 2 curves, Kani \cite{MJ} introduced the concept of a \textit{generalized Humbert scheme} $H(q)$ which is associated to a given quadratic form $q$. This set is defined by using the refined Humbert invariant as follows: Given any integral positive definite quadratic form $q$ in $r$ variables, $H(q)$ is defined by
$$
H(q) \;=\; \{\langle A,\theta\rangle  \in \mathcal{A}_2(K)\, :\,  q_{(A,\theta)} \rightarrow  q\}.
$$
The set $H(q)$ is called a \textit{generalized Humbert scheme.}   We  have that if $q(x)=Nx^2,$  then $H(q)=H_N$ is the \textit{(classical) Humbert surface of invariant $N$} (see \cite[p.~25]{MJ}). 
 Moreover, another application of this set is the determination of the components of the intersection of two Humbert surfaces $H_{N}\cap H_{M^2}$; cf.\ \cite{SubcoversofCurves}, \cite{generalizedhumbert}.

For simplicity, we use the abbreviation $ q = [a,b,c,r,s,t]$
 to denote an integral ternary quadratic form $q(x,y,z)=ax^2+by^2+cz^2+ryz+sxz+txy$, and  $q=[a,b,c]$ to denote an integral binary quadratic form $q(x,y)=ax^2+bxy+cy^2.$

Since the self-intersection number of a divisor on an abelian surface $A$ is even, i.e., $2\mid (D.D)$, for any divisor $D$ on $A$  (as in the proof of Corollary \ref{imprimitiveiffeven}), we see by definition of $ q_{(A,\theta)}$ that for any $(A,\theta)$,
\begin{equation}
    \label{eq: Humbert_surface_nonempty}
    q_{(A,\theta)}(D) \;\equiv\; 0 , 1 \Mod{4}, \ \text{ for any } D\;\in\; \NS(A,\theta).
\end{equation}
 
 We are now ready to prove Corollary \ref{cor: binaryform_H(q)} by applying our main results.


\spm
\noindent\textbf{Proof of Corollary \ref{cor: binaryform_H(q)}.}
 ($\Rightarrow$) This part follows from Eq.~(\ref{eq: Humbert_surface_nonempty}). More precisely, if $H(q)\neq \varnothing$, then there exists $(A,\theta)$ such that $q_{(A,\theta)} \rightarrow q$. Since $q_{(A,\theta)}\equiv0 , 1 \Mod{4}$ by (\ref{eq: Humbert_surface_nonempty}), it follows that also $q\equiv0,1 \Mod{4}$, which proves this implication.

 ($\Leftarrow$) We aim to construct a ternary form $f\sim q_{(A,\theta)}$, for some $(A, \theta)$ such that $f \rightarrow q$. Then this implication follows since this construction implies that $ (A,\theta) \in H(q)$. Let us put $q = [a,b,c]$. The construction of such a ternary form is done depending  on the values $a$ and $c \Mod{4}$ in three cases.
  \spm

\noindent{\bf Case 1.} Assume that $a\equiv c \equiv 0\Mod{4}$.
\sps

\noindent Let $q_1=[a_1,b_1,c_1]$ be the reduced form of $q$. Note that $H(q)=H\left(q_1\right)$ by definition since $q \sim q_1$. Let us take $f=[a_1, c_1, 4, 4,0, b_1]$. We first show that $f$ is positive definite. To see this,  note that 
\begin{equation}
    \label{eq: ternay_form_nonemptyH(q)}
    4 a_1 f \;=\; (2 a_1 x+b_1 y)^{2} \;+\; q^{\prime}(y, z)  \text {, where } q^{\prime}(y, z)\;=\;\left[4a_1c_1-b_1^{2}, 16 a_1, 16 a_1\right].
\end{equation}
Since $q_1$ is reduced, it follows from (2.12) of \cite{davidcox} that $-\disc(q_1)=4 a_1 c_1-b_1^{2} \geqslant 3 a_1^{2} \geqslant 12a_1$  (because $\left.a_1 \geqslant 4\right)$, and so $0 \geqslant 12a_1+\disc(q_1)>4 a_1+\disc(q_1)$. Thus, it follows that $\disc\left(q^{\prime}\right)=(16 a_1)^{2}+4(\disc(q_1)16a_1)=64 a_1(4 a_1+\disc(q_1))<0$. Hence, $q^{\prime}$ is positive definite (by Exercise 2.4 of \cite{davidcox}) since $q^{\prime}(1,0)=-\disc(q_1)>0$ and $\disc(q')<0$. Thus, it follows from equation (\ref{eq: ternay_form_nonemptyH(q)}) that $f$ is positive definite (since $a_1>0$).

Since $q \equiv 0,1 \Mod{4}$ (and since $a\equiv c \equiv 0 \Mod{4}$), we have that $b \equiv 0 \Mod{4}$, and so $4\mid \cont(q)=\cont(q_1)$. Thus, since $4 \mid a_1, b_1, c_1$, it follows that $\operatorname{cont}(f)=\gcd(a_1,c_1,4,4,0,b_1)=4$. Also, we see that $f/2$ is an improperly primitive form. Since $f(0,0,1)=4$, we obtain from Theorem \ref{maintheorem} that $f\sim q_{(A,\theta)}$, for some $(A,\theta)$. The assertion follows in this case because $f(x, y, 0)=a_1 x^{2}+c_1 y^{2}+b_1 x y$, i.e., $f \rightarrow q_1$.
\spm

%

\noindent{\bf Case 2.} Assume that $a\not\equiv 0 \Mod{4}$. 
\sps

\noindent Since $q\equiv 0, 1 \Mod{4}$, it follows that $a \equiv 1\Mod{4}$. Choose an odd prime number $p$ such that $p>\max \left(4 a c-b^{2}, a\right)$. Then consider $f=\left[a, c, p^{2}-a, 0,0, b\right]$. Note that $\operatorname{gcd}\left(a, p^{2}-a\right)= \operatorname{gcd}\left(a, p^{2}\right)=1$ because $p>a$. This shows that $f$ is primitive. Since $p>a$, we have that  $p^{2}-a>0$, and so  $f=a x^{2}+b x y+cy^{2}+\left(p^{2}-a\right) z^{2}$ is positive definite by the hypothesis that $q$ is positive definite.

We now claim that $\gcd(p, \disc(f)) = 1$. To see this, we first find that $\disc(f)=\left(b^{2}-4 a c\right)\left(p^{2}-a\right)$. Since $p>-b^{2}+4 a c>0$, it is clear that $p\nmid b^{2}-4 a c$. Since $\gcd(p,a)=1$, we see that $\operatorname{gcd}\left(p, p^{2}-a\right)=1$, and thus it follows that $p\nmid \disc(f)$, and the claim follows. Therefore, we obtain that $f(1,0,1)=p^{2}$ with $\operatorname{gcd}(p, \disc(f))=1$.

To conclude this case, it suffices to show that $f \equiv 0,1 \Mod{4}$. Since $p$ is odd, and $a \equiv 1 \Mod{4}, $ we have that $ p^{2}-a \equiv 0 \Mod{4}$, and thus, it follows that $f(x, y, z) \equiv a x^{2}+b x y+c y^{2} \equiv 0,1 \Mod{4}$ by the hypothesis. Therefore, the assertion follows when $a \equiv 1\Mod{4}$ by Theorem 1 of \cite{refhum} as in Case 1 since $f(x, y, 0)=q$, i.e., $f \rightarrow q$.
\sps


\noindent{\bf Case 3.} Assume that $c\not\equiv 0 \Mod{4}$. 
\sps

\noindent  We again see that $c \equiv 1\Mod{4}$ because $q\equiv 0 , 1 \Mod{4}$. Then the same argument in Case 2 applies in this case since $q \sim [c, b, a]$, and so the assertion follows. 
\spm

This corollary can be used to  prove Corollary \ref{cor: intersection_humbert_surfaces} 
as follows.
\spm

\noindent\textbf{Proof of Corollary \ref{cor: intersection_humbert_surfaces}.} By Eq.~(\ref{eq: Humbert_surface_nonempty}), it suffices to verify the second assertion  since $H_N=H(Nx^2)$ by what was mentioned above. By the hypothesis, we have that $N+M \equiv 0, 1, 2 \Mod{4}$. Suppose first that $N + M\equiv 0, 1 \Mod{4}$. Consider the binary form $q=[N, 0, M]$. It is clear that $q \equiv 0,1 \Mod{4}$. Also, it is clear that $q$ is positive definite, and thus $H(q) \neq \varnothing$ by Corollary  \ref{cor: binaryform_H(q)}. Since $q \rightarrow N, M$, we have that $\varnothing \neq H(q) \subset H_{N} \cap H_{M}$, and this proves the assertion in this case.

Assume next that $N + M\equiv 2 \Mod{4}$, which implies that  $N \equiv M \equiv 1 \Mod{4}$ by the hypothesis. We easily see from Corollary \ref{cor: binaryform_H(q)} that $H([1,0,4])\neq \varnothing$. Thus,  the statement follows when $N=M=1$, because $\varnothing\neq H([1,0,4]) \subset H_{1}=H_N \cap H_M$. So we may suppose that at least one of the numbers $N$ or $M \geqslant 5$. 

Consider the binary form $q=[N, 2, M]$. Note that $q \equiv(x+y)^{2} \Mod{4} \equiv 0,1 \Mod{4}$. Moreover, we see that it is positive definite since $q(x, y)=(x+y)^{2}+(N-1) x^{2}+(M-1) y^{2}$. Thus, $H(q)\neq \varnothing$ by Corollary \ref{cor: binaryform_H(q)}, and hence, the assertion similarly follows in this case.
\spm

The following example together with Remark \ref{rem: example_curve} below proves Proposition \ref{p: infinite_intersection_surfaces}, and discusses the properties of the element in the intersection stated in this proposition.

\begin{example}
\label{ex: infinite_intersection}
Let us choose a positive definite ternary quadratic form 
$$
q(x,y,z)= x^2 +y^2 + 3z^2 +xz + yz
$$
in Kaplansky's list; cf.\ (9) of \S5  of \cite{kaplansky}.
 This form represents \textit{all} odd positive integers, in particular, it represents all square-free odd integers. Clearly, it primitively represents all such numbers.

 Let $ \mathcal{H}$ denote the intersection
\begin{equation}
\label{ex: infinite_intersection_below}
 \bigcap_N \;\; H_{4N}, \ \text{ where } \ N \text{ runs over all positive odd square-free integers.} 
\end{equation}
We show in this example that  there exists $\langle J_C,\theta_C \rangle \in \mathcal{H}$ with the property that $q_C \sim 4q(x,y,z)$. Moreover, we show that $J_C \simeq E\times E$, where $\End(E)\simeq \Z[\sqrt{-10}]$ and $E/K$ is given by the equation $y^2 = x^3 + (140265\sqrt{5}-1023525)x + (83725920\sqrt{5}-407932200)$ and has $j$-invariant $j_E=212846400 + 95178240\sqrt{5}$.
 
To verify this, we first see by Theorem \ref{maintheorem}  that $4q$ is equivalent to a refined Humbert invariant $q_{C}$, for some curve  $C$ of genus 2 since $2q$ is improperly primitive and $4q(1,0,0)=4$.
 By using the arithmetic properties of the form $q_{C}$, we now show how to determine the geometric properties of $C$. For this, observe first that  $J_C \simeq E_1 \times E_2$, for isogenous CM elliptic curves $E_1, E_2$ by Proposition \ref{imprimitiveform} since  $q_{C}$ is a ternary form. Moreover, it follows from (\ref{determinantoftherefined})  that the discriminant of the degree map $q_{E_1, E_2}$ is $\disc(4q)/16 = -40$. We easily see (e.g.\ from \cite[p.~20]{Duncanbinaryforms}) that 
$$
q_{E_1, E_2} \ \text{ is equivalent to } \ [1,0,10] \quad \text{or} \quad [2,0,5].
$$  
We now claim that $q_{E_1, E_2} \sim [1,0,10]$ by using our results above. To prove the claim, we first find the basic invariant $I_2(q) = -160$. (Note that by Corollary \ref{imprimitiveform}, we can apply Proposition \ref{invariants}, and this also follows from this proposition without such a calculation). Since $5\mid 160$, it follows that $\chi_5$ is an assigned character of $F_q^B$ as was discussed above; (cf.\ \S\ref{ternaryforms}). We find that $F_q^B = [11, 11, 4, -4, -4, 2]$, and hence it follows that $\chi_5(F_q^B) = \chi_5(11) = 1$. Also, we can see from Proposition \ref{invariants} that $\chi_5(q_{E_1, E_2}) = \chi_5(F_q^B) = 1$. Note that $\chi_5$ is an assigned character of $q_{E_1, E_2}$ obviously since $5\mid \disc(q_{E_1, E_2}).$ This immediately proves the claim that $q_{E_1, E_2} \sim [1,0,10]$ because $\chi_5([2,0,5]) = \chi_5(2) = -1$.

Since the degree map $q_{E_1, E_2}$ represents 1, it follows that $E_1 \simeq E_2$, and thus  $J_C \simeq E \times E$, for some $E$ with $E\sim E_1$. Let us put $F = \End(E)\otimes \Q$, then by Corollary 42 of \cite{kani2011products} we find that $-40 = \disc(q_{E_1, E_2}) = \Delta_F$, where $\Delta_F$ is the discriminant of the quadratic number field $F$ and that $\End(E)$ is isomorphic to the maximal order of $F$ (since $-40$ is a fundamental discriminant). Therefore, it follows by (5.13) of \cite{davidcox} that $F\simeq \Q(\sqrt{-10})$ and $\End(E)\simeq \Z[\sqrt{-10}]$.

Since the class number $h(-40)=2$, we have two elliptic curves $E'/K$ and $E''/K$ up to isomorphism whose endomorphism ring is $\Z[\sqrt{-10}]$. It is clear that the degree maps $q_{E',E'}$ and $q_{E'',E''}$ are $\SL_2(\Z)$-equivalent to $[1,0,5]$ because there is only one binary quadratic form representing $1$ in the class group of the discriminant $-40$. This shows that $J_C\simeq E'\times E'\simeq E''\times E''$ by Theorem 67 of \cite{kani2011products}. By SageMath \cite{SageMath}, we see (by finding the roots of the Hilbert class polynomial $H_{-40}$) that the $j$-invariants of these curves are $212846400 \pm 95178240\sqrt{5}$. Since $J_C\simeq E'\times E'\simeq E''\times E''$, we may suppose that $E\simeq E'$ and the $j$-invariant $j_E$ of $E$  is $212846400 + 95178240\sqrt{5}$. By SageMath \cite{SageMath}, we give a (minimal) model for this curve $E$, and obtain that $E$ is given by an equation $y^2 = x^3 + (140265\sqrt{5}-1023525)x + (83725920\sqrt{5}-407932200)$,  which proves the assertion above.
\end{example}

 In Example \ref{ex: infinite_intersection}, we constructed an element in the intersection $\mathcal{H}$, and determined the geometric properties of this element. We can actually say more about this element in the following remark, and so Proposition \ref{p: infinite_intersection_surfaces} follows.

\begin{remark} 
\label{rem: example_curve} 
 Let $C$ be a curve such that $q_{(J_C,\theta_C)}\sim 4q$, where $q$ is as in Example \ref{ex: infinite_intersection}. We now show that $C$ is isomorphic to the genus $2$ curve given by the polynomial equation
 \begin{equation}
 \label{eq: poly_eq_curve}
     y^2 \;=\; 9x^5 \;+\; 40x^3 \;+\; 45x.
 \end{equation}
 For this, we use Hayashida's formula (see \S7 of \cite{hayashida1968class}) and the table of \cite[p.~271]{gelin2019principally}.

 Since $J_C\simeq E\times E$, where $\End(E)$ is a maximal order (see Example \ref{ex: infinite_intersection}), one can apply Hayashida's formula (see \S7 of \cite{hayashida1968class}) to find the number of representatives of the principal polarizations on $J_C$, and so it follows from this formula that there are two such representatives. 

 Consider the basis $\{1,\sqrt{-10}\}$ of $\Z[\sqrt{-10}]$, and let $\{1_E,\alpha\}$ be the corresponding basis of $\End(E)\simeq \Z[\sqrt{-10}]$. Then the degree map is $q(x+y\sqrt{-10})=[1,0,10]$ (see \S3.4 of \cite{kani2011products}) under this basis. Thus $h=1_E+\alpha$ is an isogeny of degree $11$.

We see by Corollary 25 of \cite{MJ} that  $\theta:= \textbf{D}(2,6,h)$ and $\theta':=\textbf{D}(3,4,h)$ are in $\mathcal{P}(J_C)$. Moreover, it follows from Corollary \ref{eventhetaprop} that $\theta \in \mathcal{P}(J_C)^{\ev}$ (since $4 \mid \disc(q_{E_1, E_2}) = -40$) and that $\theta' \in \mathcal{P}(J_C)^{\odd}$. Thus, we see from Corollary \ref{imprimitiveiffeven} that $q_{(J_C,\theta)}$ is imprimitive and $q_{(J_C,\theta')}$ is primitive. In particular, they are not equivalent forms, which means that $\theta$ and $\theta'$ are not isomorphic. Since $\mathcal{P}(J_C)$ has only two representatives by what was mentioned above, we conclude that all their representatives are $\theta$ and $\theta'$. Thus, $\mathcal{P}(J_C)^{\ev}$ has a unique representative, namely $\theta$. Since $q_{(J_C,\theta_C)}$ is imprimitive, it follows from Corollary \ref{imprimitiveiffeven} that $\theta_C\in \mathcal{P}(J_C)^{\ev}$, and so, it follows that $\theta_C\simeq \theta$.

 We know by Theorem 5.2.4 of \cite{lange2013complex} that the principal polarization $\theta_1=\textbf{D}(2,6,1_E+\alpha)$ corresponds to a positive definite unimodular Hermitian matrix with the entries in $\Z[\sqrt{-10}]$. By Proposition 61 of \cite{MJ}, we see that this matrix is 
     $\big(\begin{smallmatrix}
 2 & 1+\sqrt{-10} \\
1-\sqrt{-10} & 6
\end{smallmatrix}\big)$  since $1_E+\alpha$ corresponds to $1+\sqrt{-10}$. Note that this matrix is the one listed in Table 4 of \cite{gelin2019principally}, and it follows from this table that $\theta_1 \simeq \theta_C$ is isomorphic to a genus 2 curve given by the polynomial equation in (\ref{eq: poly_eq_curve}), as claimed.

 As a consequence, the principally polarized abelian surface $(J_C, \theta_C)$ corresponding to the curve $C$ given by $y^2 \;=\; 9x^5 \;+\; 40x^3 \;+\; 45x$ lies in the intersection $\mathcal{H}$ as claimed in Proposition \ref{p: infinite_intersection_surfaces}.
\end{remark}

\begin{remark}
It may be useful to note a generalization of Example \ref{ex: infinite_intersection}.   Kaplansky \cite{kaplansky} proved that there are at most $23$ positive definite ternary quadratic forms representing all positive odd integers. He listed all such forms and proved that $19$ of the $23$ represent all odd positive integers, and these are the forms having line numbers $1\text{\textendash}19$ in his paper. Then he remarked that the remaining four forms in his list are “\textit{plausible candidates}” to represent all positive odd integers, so these are the forms in the lines from $20\text{\textendash}23$. Then, Jagy \cite{Jagy} proved that one of Kaplansky's candidates, namely line $22$, represents all positive odd integers. Furthermore, Rouse has proved (see Theorem 7 of \cite{Jeremy_Rouse}) that, assuming GRH (the Generalized Riemann Hypothesis), the remaining three Kaplansky candidates represent all positive odd integers.  Therefore, we have a complete list of positive definite ternary quadratic forms representing all positive odd integers under GRH.

For any form $q$ in Kaplansky's list (excluding those in 1-5), it is evident that $4q \sim q_{C}$ for some curve $C$ of genus 2 by Theorem \ref{maintheorem}. Thus it follows that $(J_C, \theta_C) \in \mathcal{H}$, for all $C$ such that $ q_{C} \sim 4q$, where $q$ is a form in the list except in 20, 21, and 23 (and 1-5) without assuming GRH. If we assume GRH, then  we also have that $(J_C, \theta_C) \in \mathcal{H}$, for all $C$ such that $ q_{C} \sim 4q$, where $q$ is in 20, 21 and 23 of the list.
\end{remark}

It would be very interesting to prove (or disprove) geometrically that we have $q_{C} \sim 4f $, for some $(J_C, \theta_C) \in \mathcal{H}$, for the form $f$ in the Kaplansky's list numbered 20, 21 or 23 without assuming GRH.
\spm

\section*{Acknowledgments}
 I express my deepest gratitude to my supervisor, Prof.\ Ernst Kani, for all his help and support. I especially want to thank him for two reasons. First, he shared his current (preprint) papers with me, and he also provided me with some useful results as indicated each time. Second, he carefully read all the paper, and provided me with a long list of errors and suggestions. I also want to thank Eda K{\i}r{\i}ml{\i} for proofreading my manuscript numerous times and giving many suggestions on my writing, and  pointing out my errors and typos. I finally would like to thank the referee for his/her valuable comments, suggestions, and lists of corrections.

\bibliographystyle{plainurl}
\bibliography{thesis}

@incollection {MitaniShioda,
    AUTHOR = {Shioda, T. and Mitani, N.},
     TITLE = {Singular abelian surfaces and binary quadratic forms},
 BOOKTITLE = {Classification of algebraic varieties and compact complex
              manifolds},
     PAGES = {259--287. Lect. Notes in Math., Vol. 412},
     PUBLISHER = {Springer, Berlin},
     YEAR = {(1974)},
   MRCLASS = {14J20 (14K15)},
  MRNUMBER = {0382289},
MRREVIEWER = {P. Roquette},
URL = {https://doi.org/10.1007/BFb0066163
},
}

@book {Duncanbinaryforms,
    AUTHOR = {Buell, D. A.},
     TITLE = {Binary Quadratic Forms. Classical Theory and Modern Computations},
      
 PUBLISHER = {Springer-Verlag, New York},
      YEAR = {(1989)},
    PAGES = {x+247},
      OPTISBN = {0-387-97037-1},
   MRCLASS = {11E16 (11-02)},
  MRNUMBER = {1012948},
MRREVIEWER = {A. G. Earnest},
       OPTDOI = {10.1007/978-1-4612-4542-1},
       URL = {https://doi.org/10.1007/978-1-4612-4542-1},
}

@article {brandt1951zahlentheorie,
    AUTHOR = {Brandt, H.},
     TITLE = {Zur {Z}ahlentheorie der tern\"{a}ren quadratischen {F}ormen},
   JOURNAL = {Math. Ann.},
  FJOURNAL = {Mathematische Annalen},
    VOLUME = {124},
      YEAR = {(1952)},
     PAGES = {334--342},
      OPTISSN = {0025-5831,1432-1807},
   MRCLASS = {10.0X},
  MRNUMBER = {51269},
MRREVIEWER = {R.\ Hull},
       OPTDOI = {10.1007/BF01343574},
       URL = {https://doi.org/10.1007/BF01343574},
}

@article {brandt1952mass,
    AUTHOR = {Brandt, H.},
     TITLE = {\"{U}ber das {M}ass positiver tern\"{a}rer quadratischer
              {F}ormen},
   JOURNAL = {Math. Nachr.},
  FJOURNAL = {Mathematische Nachrichten},
    VOLUME = {6},
      YEAR = {(1952)},
     PAGES = {315--318},
      OPTISSN = {0025-584X,1522-2616},
   MRCLASS = {10.0X},
  MRNUMBER = {51268},
MRREVIEWER = {R.\ Hull},
       OPTDOI = {10.1002/mana.19520060507},
       URL = {https://doi.org/10.1002/mana.19520060507},
}

@article {kani2011products,
    AUTHOR = {Kani, E.},
     TITLE = {Products of {CM} elliptic curves},
   JOURNAL = {Collect. Math.},
  FJOURNAL = {Collectanea Mathematica},
    VOLUME = {62},
      YEAR = {(2011)},
    NUMBER = {3},
     PAGES = {297--339},
      OPTISSN = {0010-0757},
   MRCLASS = {11G10 (11G15 14H52 14K02 14K22 14L15)},
  MRNUMBER = {2825715},
MRREVIEWER = {Joseph H. Silverman},
       OPTDOI = {10.1007/s13348-010-0029-1},
        URL = {https://doi.org/10.1007/s13348-010-0029-1},
}

@article {kani2014jacobians,
    AUTHOR = {Kani, E.},
     TITLE = {Jacobians isomorphic to a product of two elliptic curves and ternary quadratic forms},
   JOURNAL = {J. Number Theory},
  FJOURNAL = {Journal of Number Theory},
    VOLUME = {139},
      YEAR = {(2014)},
     PAGES = {138--174},
      OPTISSN = {0022-314X,1096-1658},
   MRCLASS = {14H40 (11G10 11G15 11G18 14H30 14H52)},
  MRNUMBER = {3173190},
MRREVIEWER = {Edward\ F.\ Schaefer},
       OPTDOI = {10.1016/j.jnt.2013.12.006},
       URL = {https://doi.org/10.1016/j.jnt.2013.12.006},
}

@article {ESCI,
    AUTHOR = {Kani, E.},
     TITLE = {Elliptic subcovers of a curve of genus 2. {I}. {T}he isogeny defect},
   JOURNAL = {Ann. Math. Qu\'{e}.},
  FJOURNAL = {Annales Math\'{e}matiques du Qu\'{e}bec},
    VOLUME = {43},
      YEAR = {(2019)},
    NUMBER = {2},
     PAGES = {281--303},
      OPTISSN = {2195-4755,2195-4763},
   MRCLASS = {14H30 (14H05 14H25 14H40)},
  MRNUMBER = {3996071},
MRREVIEWER = {Giancarlo\ Urz\'{u}a},
       OPTDOI = {10.1007/s40316-018-0105-6},
       URL = {https://doi.org/10.1007/s40316-018-0105-6},
}

@article {ESCII,
    AUTHOR = {Kani, E.},
     TITLE = {Elliptic subcovers of a curve of genus 2 {II}. {T}he refined
              {H}umbert invariant},
   JOURNAL = {J. Number Theory},
  FJOURNAL = {Journal of Number Theory},
    VOLUME = {193},
      YEAR = {(2018)},
     PAGES = {302--335},
      OPTISSN = {0022-314X,1096-1658},
   MRCLASS = {14H30 (11G30 14H40)},
  MRNUMBER = {3846811},
MRREVIEWER = {Sajad\ Salami},
       OPTDOI = {10.1016/j.jnt.2018.05.011},
       URL = {https://doi.org/10.1016/j.jnt.2018.05.011},
}

@incollection {SubcoversofCurves,
    AUTHOR = {Kani, E.},
     TITLE = {Subcovers of curves and moduli spaces},
 BOOKTITLE = {Geometry at the frontier},
    SERIES = {Contemp. Math.},
    VOLUME = {766},
     PAGES = {229--250},
 PUBLISHER = {Amer. Math. Soc., RI},
      YEAR = {(2021)},
      OPTISBN = {978-1-4704-5327-5},
   MRCLASS = {14G35},
  MRNUMBER = {4248056},
MRREVIEWER = {Nicolae\ Manolache},
       OPTDOI = {10.1090/conm/766/15384},
       URL = {https://doi.org/10.1090/conm/766/15384},
}

@article{generalizedhumbert,
  title={Generalized {Humbert} schemes and intersections of {Humbert} surfaces},
  author={Kani, E.},
  year={2021},
  note={Preprint, 35 pages},
  URL={https://mast.queensu.ca/~kani/papers/interHum11.pdf}
}

@article{refhum,
  title={The refined {Humbert} invariant for abelian product surfaces with complex multiplication},
  author={Kani, E.},
  year={2025},
  note={Preprint, 24 pages},
  URL={https://mast.queensu.ca/~kani/papers/refHum5.pdf}
}

@article {MJ,
    AUTHOR = {Kani, E.},
     TITLE = {The moduli spaces of {J}acobians isomorphic to a product of
              two elliptic curves},
   JOURNAL = {Collect. Math.},
  FJOURNAL = {Collectanea Mathematica},
    VOLUME = {67},
      YEAR = {(2016)},
     PAGES = {21--54},
      OPTISSN = {0010-0757},
   MRCLASS = {14H10 (14H40)},
  MRNUMBER = {3439838},
MRREVIEWER = {Francisco J. Plaza Mart\'{\i}n},
       OPTDOI = {10.1007/s13348-015-0148-9},
       URL = {https://doi.org/10.1007/s13348-015-0148-9},
}

@article {Mcmullen,
    AUTHOR = { C. McMullen},
     TITLE = {Teichm\"{u}ller curves in genus two: discriminant and spin},
   JOURNAL = {Math. Ann.},
  FJOURNAL = {Mathematische Annalen},
    VOLUME = {333},
      YEAR = {(2005)},
    NUMBER = {1},
     PAGES = {87--130},
      OPTISSN = {0025-5831,1432-1807},
   MRCLASS = {32G15 (30F10 37D50)},
  MRNUMBER = {2169830},
MRREVIEWER = {Thomas\ A.\ Schmidt},
      OPTDOI = {10.1007/s00208-005-0666-y},
       URL = {https://doi.org/10.1007/s00208-005-0666-y},
}

@book {van2012hilbert,
    AUTHOR = {van der Geer, G.},
     TITLE = {Hilbert modular surfaces},
    SERIES = {\emph{Ergebnisse der Mathematik und ihrer Grenzgebiete (3)}},
    VOLUME = {16},
 PUBLISHER = {Springer-Verlag, Berlin},
      YEAR = {(1988)},
     PAGES = {x+291},
      OPTISBN = {3-540-17601-2},
   MRCLASS = {11F41 (11G10 11G15 14J20)},
  MRNUMBER = {930101},
MRREVIEWER = {O.\ V.\ Shvartsman},
       OPTDOI = {10.1007/978-3-642-61553-5},
       URL = {https://doi.org/10.1007/978-3-642-61553-5},
}

@article{hayashida1968class,
  title={A class number associated with the product of an elliptic curve with itself},
  author={Hayashida, T.},
  journal={Journal of the Mathematical Society of Japan},
  volume={20},
  pages={26--43},
  year={(1968)},
  publisher={The Mathematical Society of Japan}
}

@inproceedings {gelin2019principally,
    AUTHOR = {G\'{e}lin, A. and Howe, E. and Ritzenthaler,
              C.},
     TITLE = {Principally polarized squares of elliptic curves with field of
              moduli equal to {$\mathbb{Q}$}},
 BOOKTITLE = {Proc. 13th {A}lgorithmic {N}umber
              {T}heory {S}ymposium},
    SERIES = {Open Book Ser. 2},
    VOLUME = {},
     PAGES = {257--274},
 PUBLISHER = {},
      YEAR = {(2019)},
      OPTISBN = {978-1-935107-03-3; 978-1-935107-02-6},
   MRCLASS = {11G15 (14H25 14H40)},
  MRNUMBER = {3952016},
MRREVIEWER = {James\ H.\ Stankewicz},
URL = {https://msp.org/obs/2019/2-1/p16.xhtml},
}

@article {kani1994elliptic,
    AUTHOR = {Kani, E.},
     TITLE = {Elliptic curves on abelian surfaces},
   JOURNAL = {Manuscripta Math.},
  FJOURNAL = {Manuscripta Mathematica},
    VOLUME = {84},
      YEAR = {(1994)},
    NUMBER = {2},
     PAGES = {199--223},
      OPTISSN = {0025-2611,1432-1785},
   MRCLASS = {14K10 (14J25 14K20)},
  MRNUMBER = {1285957},
MRREVIEWER = {C.\ A. M. Peters},
       OPTDOI = {10.1007/BF02567454},
       URL = {https://doi.org/10.1007/BF02567454},
}

@book {watson1960integral,
    AUTHOR = {Watson, G. L.},
     TITLE = {Integral quadratic forms},
    OPTSERIES = {Cambridge Tracts in Mathematics and Mathematical Physics},
    OPTVOLUME = {No. 51},
 PUBLISHER = {Cambridge U. Press, Cambridge},
      YEAR = {(1960)},
     PAGES = {xii+143},
   MRCLASS = {10.00},
  MRNUMBER = {118704},
MRREVIEWER = {B.\ W.\ Jones},
}

@book{dicksonsbook,
  title={Studies in Number Theory},
  author={ Dickson, L.},
  year={(1957)},
  publisher={U Chicago Press, Chicago, 1930.
Reprinted by Chelsea Publ. Co., New York}
}

@book{smith,
  title={On the orders and genera of ternary quadratic forms \emph{{(1867)}}},
  author={Smith, J. H. S.},
  year={(1894)},
  publisher={In: Collect. Math. Papers vol. I, Oxford, pp. 455\textsc{\textendash}509}
}

@book{jones1950arithmetic,
  title={The Arithmetic theory of Quadratic Forms},
  author={Jones, B.},
  year={(1950)},
  publisher={Carus Math. Monogr., Wiley, New York},
  URL ={https://www.jstor.org/stable/10.4169/j.ctt5hh98f}
}

@book {mumford1970abelian,
    AUTHOR = {Mumford, D.},
     TITLE = {Abelian varieties},
    OPTSERIES = {Tata Institute of Fundamental Research Studies in Mathematics},
    OPTVOLUME = {5},
 PUBLISHER = {2nd edn. Oxford University Press, London},
      YEAR = {(1970)},
     PAGES = {viii+242},
   MRCLASS = {14.51},
  MRNUMBER = {282985},
MRREVIEWER = {James\ Milne},
}

@book{lange2013complex,
  title={Complex abelian varieties},
  author={Lange, H. and Birkenhake, C.},
  volume={302},
  year={2013},
  publisher={Springer Science}
}

@incollection {milne1986jacobian,
    AUTHOR = {Milne, J. S.},
     TITLE = {Jacobian varieties},
 BOOKTITLE = {Arithmetic geometry ({S}torrs, {C}onn., 1984)},
     PAGES = {167--212},
 PUBLISHER = {Springer, New York},
      YEAR = {1986},
      ISBN = {0-387-96311-1},
   MRCLASS = {14H40},
  MRNUMBER = {861976},
}

@book {davidcox,
    AUTHOR = {Cox, D. A.},
     TITLE = {Primes of the form {$x^2 + ny^2$}},
    SERIES = {A Wiley-Interscience Publ.},
      OPTNOTE = {Fermat, class field theory and complex multiplication},
 PUBLISHER = {John Wiley \& Sons, Inc., New York},
      YEAR = {(1989)},
     PAGES = {xiv+351},
      OPTISBN = {0-471-50654-0; 0-471-19079-9},
   MRCLASS = {11A41 (11F11 11R11 11R16 11R18 11R37 11Y11)},
  MRNUMBER = {1028322},
MRREVIEWER = {Andrew\ Bremner},
URL = {https://onlinelibrary.wiley.com/doi/book/10.1002/9781118032756},
}

@Manual{SageMath,
  key          = {SageMath},
  author       = {The Sage Developers},
  title        = {SageMath, the Sage Mathematics Software System (Version 9.6)},
  note         = {\url{https://www.sagemath.org}},
  year         = {2022},
}

@article {auffarth,
    AUTHOR = {Auffarth, II, R.},
     TITLE = {Elliptic curves on {A}belian varieties},
   JOURNAL = {Illinois J. Math.},
  FJOURNAL = {Illinois Journal of Mathematics},
    VOLUME = {59},
      YEAR = {(2015)},
    NUMBER = {2},
     PAGES = {319--336},
      OPTISSNN = {0019-2082,1945-6581},
   MRCLASS = {14K12 (14H52 14K02 32G20)},
  MRNUMBER = {3499514},
MRREVIEWER = {Giulio\ Codogni},
       URL = {http://projecteuclid.org/euclid.ijm/1462450703},
}

@article {Runge,
    AUTHOR = {Runge, B.},
     TITLE = {Endomorphism rings of abelian surfaces and projective models
              of their moduli spaces},
   JOURNAL = {Tohoku Math. J. (2)},
  FJOURNAL = {The Tohoku Mathematical Journal. Second Series},
    VOLUME = {51},
      YEAR = {(1999)},
    NUMBER = {3},
     PAGES = {283--303},
      OPTISSN = {0040-8735,2186-585X},
   MRCLASS = {14K10 (14G35 14K22)},
  MRNUMBER = {1707758},
MRREVIEWER = {Bruce\ Hunt},
       OPTDOI = {10.2748/tmj/1178224764},
       URL = {https://doi.org/10.2748/tmj/1178224764},
}

@article {quaternionic,
    AUTHOR = {Lin, Y. and Yang, Y.},
     TITLE = {Quaternionic loci in {S}iegel's modular threefold},
   JOURNAL = {Math. Z.},
  FJOURNAL = {Mathematische Zeitschrift},
    VOLUME = {295},
      YEAR = {(2020)},
    NUMBER = {1-2},
     PAGES = {775--819},
      OPTISSN = {0025-5874,1432-1823},
   MRCLASS = {11G15 (11F03 11F46 11G10)},
  MRNUMBER = {4100031},
MRREVIEWER = {Nathan\ Grieve},
       OPTDOI = {10.1007/s00209-019-02372-z},
       URL = {https://doi.org/10.1007/s00209-019-02372-z},
}

@article {Hashimoto,
    AUTHOR = {Hashimoto, K.},
     TITLE = {Explicit form of quaternion modular embeddings},
   JOURNAL = {Osaka J. Math.},
  FJOURNAL = {Osaka Journal of Mathematics},
    VOLUME = {32},
      YEAR = {(1995)},
    NUMBER = {3},
     PAGES = {533--546},
      OPTISSN = {0030-6126},
   MRCLASS = {11F46 (11F32)},
  MRNUMBER = {1367889},
MRREVIEWER = {Jae-Hyun\ Yang},
       URL = {http://projecteuclid.org/euclid.ojm/1200786264},
}

@article {Pardini_Marco_Rollenske,
    AUTHOR = {Franciosi, M. and Pardini, R. and Rollenske, S.},
     TITLE = {{$(d,d')$}-elliptic curves of genus two},
   JOURNAL = {Ark. Mat.},
  FJOURNAL = {Arkiv f\"{o}r Matematik},
    VOLUME = {56},
      YEAR = {(2018)},
    NUMBER = {2},
     PAGES = {299--317},
      OPTISSN = {0004-2080,1871-2487},
   MRCLASS = {14H30 (14H52)},
  MRNUMBER = {3893776},
MRREVIEWER = {Valentijn\ Zo\"{e}\ Karemaker},
       OPTDOI = {10.4310/ARKIV.2018.v56.n2.a6},
       URL = {https://doi.org/10.4310/ARKIV.2018.v56.n2.a6},
}

@article {kaplansky,
    AUTHOR = {Kaplansky, I.},
     TITLE = {Ternary positive quadratic forms that represent all odd
              positive integers},
   JOURNAL = {Acta Arith.},
  FJOURNAL = {Acta Arithmetica},
    VOLUME = {70},
      YEAR = {(1995)},
    NUMBER = {3},
     PAGES = {209--214},
      OPTISSN = {0065-1036,1730-6264},
   MRCLASS = {11E25 (11D85)},
  MRNUMBER = {1322563},
MRREVIEWER = {Jorge\ F.\ Morales},
      OPTDOI = {10.4064/aa-70-3-209-214},
       URL = {https://doi.org/10.4064/aa-70-3-209-214},
}

@article {Jeremy_Rouse,
    AUTHOR = {Rouse, J.},
     TITLE = {Quadratic forms representing all odd positive integers},
   JOURNAL = {Amer. J. Math.},
  FJOURNAL = {American Journal of Mathematics},
    VOLUME = {136},
      YEAR = {(2014)},
    NUMBER = {6},
     PAGES = {1693--1745},
      OPTISSN = {0002-9327,1080-6377},
   MRCLASS = {11E20 (11D85 11E12 11E25 11H55)},
  MRNUMBER = {3282985},
MRREVIEWER = {Laurent\ Habsieger},
       OPTDOI = {10.1353/ajm.2014.0041},
       URL = {https://doi.org/10.1353/ajm.2014.0041},
}

@article {Jagy,
    AUTHOR = {Jagy, W. C.},
     TITLE = {Five regular or nearly-regular ternary quadratic forms},
   JOURNAL = {Acta Arith.},
  FJOURNAL = {Acta Arithmetica},
    VOLUME = {77},
      YEAR = {1996},
    NUMBER = {4},
     PAGES = {361--367},
      OPTISSN = {0065-1036,1730-6264},
   MRCLASS = {11E25 (11E20)},
  MRNUMBER = {1414516},
MRREVIEWER = {A.\ G.\ Earnest},
       OPTDOI = {10.4064/aa-77-4-361-367},
       URL = {https://doi.org/10.4064/aa-77-4-361-367},
}



\end{document}